\documentclass[11pt]{amsart}
\usepackage{amsxtra}

\theoremstyle{plain}

\def\CC {{\mathbb C}}
\def\RR {{\mathbb R}}
\def\NN {{\mathbb N}}
\def\ZZ {{\mathbb Z}}
\def\PP {{\mathbb P}}

\def\be {\begin{eqnarray}}
\def\ben {\begin{eqnarray*}}
\def\ee {\end{eqnarray}}
\def\een {\end{eqnarray*}}

\def\AAA{\kern-0.3em}
\def\AA{\kern-0.18em}
\def\AC{\kern-0.14em}
\def\AB{\kern-0.22em}

\newcommand \nc {\newcommand}

\newtheorem{theorem}{Theorem}[section]
\newtheorem{lemma}[theorem]{Lemma}
\newtheorem{proposition}[theorem]{Proposition}
\newtheorem{corollary}[theorem]{Corollary}
\newtheorem{definition}[theorem]{Definition}
\newtheorem{example}[theorem]{Example}
\newtheorem{remark}[theorem]{Remark}
\newtheorem{conjecture}[theorem]{Conjecture}
\nc \bth[1] { \begin{theorem}\label{t#1} } \nc \ble[1] {
\begin{lemma}\label{l#1} } \nc \bpr[1] {
\begin{proposition}\label{p#1} } \nc \bco[1] {
\begin{corollary}\label{c#1} } \nc \bde[1] {
\begin{definition}\label{d#1}\rm } \nc \bex[1] {
\begin{example}\label{e#1}\rm } \nc \bre[1] {
\begin{remark}\label{r#1}\rm } \nc \bcon[1] { 
\begin{conjecture}\label{co#1}\rm} \nc \bqu[1]  {
\medskip\noindent{\it{Question #1}} }
\nc {\ethe} { \end{theorem} }
 \nc {\ele} { \end{lemma} } \nc {\epr}
{ \end{proposition} } \nc {\eco} { \end{corollary} } \nc {\ede} {
\end{definition} } \nc {\eex} { \end{example} } \nc {\ere} {
\end{remark} } \nc {\econ} { \end{conjecture} } \nc {\equ} {\smallskip}
 \nc \thref[1]{Theorem \ref{t#1}}
\nc \leref[1]{Lemma \ref{l#1}} \nc \prref[1]{Proposition
\ref{p#1}} \nc \coref[1]{Corollary \ref{c#1}} \nc
\deref[1]{Definition \ref{d#1}} \nc \exref[1]{Example \ref{e#1}}
\nc \reref[1]{Remark \ref{r#1}}
\nc \conref[1]{Conjecture \ref{co#1}}

\def \B {{\mathcal B}}

\def \L {{\mathcal L}}

\def \diag { {\mathrm{diag}} }
\def \res { {\mathrm{Res}} }

\def \GL { {\mathrm{GL}}}

 \def\AA  {\kern-0.1em}
 \def\BB  {\kern+0.1em}
 \def\BBB {\kern+0.15em}
 \def\K   {\kern+0.05em}
 \def\MK  {\kern-0.07em}
 \def\MKK {\kern-0.04em}

\begin{document}

\vspace{0.5cm}

\title[Zero  level perturbation]
{ Zero  level perturbation of a certain third-order linear solvable
  ODE with an irregular singularity \\at the origin of Poincar\'e rank 1 }

\author[Tsvetana  Stoyanova]{Tsvetana  Stoyanova}

\date{06.05.2018}

 \maketitle

\begin{center}
{Department of Mathematics and Informatics,
Sofia University,\\ 5 J. Bourchier Blvd., Sofia 1164, Bulgaria, 
cveti@fmi.uni-sofia.bg}
\end{center}

\vspace{0.5cm}

{\bf Abstract.}  We study  an irregular 
singularity of Poincar\'e rank 1  at the origin of a certain third-order
linear solvable homogeneous ODE. We perturb the equation by introducing  a small 
parameter $\varepsilon\in(\RR_+,0)$ ($\varepsilon < 1$), which causes the splitting of the irregular singularity
into two finite Fuchsian singularities.   
 We show that when the solutions of the perturbed equation contain logarithmic
terms, the Stokes matrices  of the initial equation are limits
of the part of  the monodromy matrices around the finite resonant Fuchsian singularities
of the perturbed equation.\\



\headsep 10mm \oddsidemargin 0in \evensidemargin 0in

\section{Introduction}

    We consider a linear complex ordinary differential equation  
	  \be\label{initial}
		   L\,y=0\,,
			\ee
			where $L$ is a third-order linear differential operator of the form
			\be\label{form-0}
			  L=L_3 \circ L_2 \circ L_1
			\ee
			with
			\be\label{npe}
		  L_1 =\partial-\frac{1}{x^2}\,,\quad
			L_2 =\partial - \frac{\nu-2}{x} - \frac{2}{x^2}\,,\quad
			L_3 =\partial - \frac{\nu-4}{x}\,,\quad \partial=\frac{d}{d x}
		\ee
		and $\nu\in\CC$.
			The equation \eqref{initial} with \eqref{form-0}-\eqref{npe} (in short, the initial equation) 
			is a third-order solvable differential equation, in the sense that its differential Galois
		group is a solvable linear algebraic group \cite{K}. We consider the equation \eqref{initial} 
		over $\CC\PP^1$, where it has an irregular singular point at the origin
	  of Poincar\'e rank 1 and a regular singular point (if $\nu \neq 0$, see \reref{E1})
		at $x=\infty$. We associate with the initial equation its analytic invariants at
		the irregular singularity - the Stokes matrices and the formal monodromy.

			On the other hand, 
		following \cite{AG,AG1, HLR, MK,CL-CR,CL-CR1}, we can always regard the irregular
		singularity at the origin as a result of confluence of two Fuchsian singularities.
		Namely, introducing a (small) parameter $\varepsilon\in (\RR_+, 0)$,
			we consider a perturbation of the initial equation  
			\be\label{perturbed}
			 L(\varepsilon)\,y=0\,,
		  \ee
		where $L(\varepsilon)$ is again a third-order differential operator of the form
		\be\label{form}
			L(\varepsilon)=L_3(\varepsilon) \circ L_2(\varepsilon) \circ
			                  L_1(\varepsilon)\,.						
		\ee
		The first-order differential operators $L_j(\varepsilon),\,j=1, 2, 3$
		are defined as follows, 
		  \be\label{pe}
			 & &
		    L_1(\varepsilon)=\partial - \frac{1}{2 \sqrt{\varepsilon}}\,
				\left(\frac{1}{x-\sqrt{\varepsilon}}
					- \frac{1}{x+\sqrt{\varepsilon}}\right)\,,\nonumber\\[0.9ex]
				& &
				L_2(\varepsilon)=\partial - \left(\frac{\nu-2}{2} + \frac{1}{\sqrt{\varepsilon}}\right)\,
					               \frac{1}{x-\sqrt{\varepsilon}} -
				\left(\frac{\nu-2}{2} - \frac{1}{\sqrt{\varepsilon}}\right)\,
												\frac{1}{x+\sqrt{\varepsilon}}\,,\\[0.9ex]
				& &									
			L_3(\varepsilon)=\partial - \frac{\nu-4}{2}\,
				\left(\frac{1}{x-\sqrt{\varepsilon}} + \frac{1}{x+\sqrt{\varepsilon}}\right)\,,\nonumber\\[0.9ex]
				& &
			L_j(0)=L_j\,,j=1, 2, 3\,.\nonumber	
							  \ee
			 The equation \eqref{perturbed} with \eqref{form}-\eqref{pe} 
			(in short, the perturbed equation) is a third-order Fuchsian equation.					
				It has over $\CC\PP^1$ three Fuchsian singularities: $x=-\sqrt{\varepsilon},\,
				\sqrt{\varepsilon},\,\infty$. Following \cite{CL-CR, CL-CR1},
			through this article, we denote the finite regular singularities by
				$x_L=-\sqrt{\varepsilon}$ and $x_R=\sqrt{\varepsilon}$ for $\sqrt{\varepsilon}\in\RR_+$.
      In this article we consider the perturbed equation together with its
			monodromy matrices around the finite Fuchsian singularities $x_L$ and
			$x_R$.
			
			The main result of this article is, that in the presence of logarithmic terms
			in the solutions of the perturbed equation, the Stokes matrices  of the initial equation
			are limits of the parts of the monodromy matrices around the
			resonant singular points $x_j, j=R, L$ 	of the perturbed equation.
      More precisely, it turns out that exactly these parts of the monodromy matrices,
			which govern the presence of logarithmic terms tend to the Stokes matrices
			when $\sqrt{\varepsilon}$ tends to zero.

					The point $x=\infty$ is a regular singularity for both initial and perturbed equations
					(except the case when $\nu=0$, see \reref{E1} and \reref{E2}).
					Furthermore the exponents $\rho^{\infty}_1=0,\,\rho^{\infty}_2=1-\nu,\,
					\rho^{\infty}_3=2-\nu$ at $x=\infty$ are the same for the both equations.
		      The perturbation, defined above, splits the irregular singularity
		$x=0$ of Poincar\'e rank 1 of the initial  equation into two finite Fuchsian singular points
		$x_R=\sqrt{\varepsilon}$ and $x_L=-\sqrt{\varepsilon}$ of the perturbed equation, but does not change
		 the infinity point. 
		There exists a more generic perturbation, which not only
		splits the irregular singularity, but also perturbs the exponents 
		$\rho^{\infty}_j,\,j=1, 2, 3$ at $x=\infty$. 					
		   Generally, if the  coefficients $a_j(x)\,,\,j=1, 2, 3$ of
		the initial equation  \eqref{form-g} are given by
		 $$\,
		    	a_j(x)=\frac{\alpha_j}{x} + \frac{\beta_j}{x^2}\,,\quad
					\alpha_j,\,\beta_j\in\CC\quad
					\textrm{such that }\quad (\beta_1, \beta_2, \beta_3) \neq (0, 0, 0)\,,
		\,$$
		 then the corresponding perturbed coefficients
		$a_j(x, \varepsilon)$ would be 
      \be\label{pco}
		   a_j(x, \varepsilon)
			    &=&
			\left(\frac{\alpha_j(\varepsilon)}{2} + \frac{\beta_j(\varepsilon)}{2 \sqrt{\varepsilon}}\right)\,
			     \frac{1}{x-\sqrt{\varepsilon}} + 
			\left(\frac{\alpha_j(\varepsilon)}{2} - \frac{\beta_j(\varepsilon)}{2 \sqrt{\varepsilon}}\right)\,		
			\frac{1}{x+\sqrt{\varepsilon}}\,,\\[0.9ex]
			  a_j(x, 0) &=& a_j(x)\,,\nonumber
		\ee
		where $\alpha_j(\varepsilon)$ and $\beta_j(\varepsilon)$ are polynomials in
		$\varepsilon$ such that $\alpha_j(0)=\alpha_j$ and $\beta_j(0)=\beta_j$.
		We call such a perturbation a $N$-th level perturbation if all these polynomials
		are of degree $N$. Then the $N$-th  level perturbation has the exponents
	 $\rho^{\infty}_1=-A_1(\varepsilon),\,
		\rho^{\infty}_2=1-\nu-A_2(\varepsilon),\,\rho^{\infty}_3=2-\nu-A_3(\varepsilon)$,
		where $A_j(\varepsilon)$ are polynomials in $\varepsilon$ of degree $N$ such that
		$A_j(0)=0$ and $\alpha_j(\varepsilon)=\alpha_j+A_j(\varepsilon)$.
		So the zero level perturbation is a perturbation
		with constant polynomials $\alpha_j$ and $\beta_j$.

     Recently several authors (see below) have studied the irregular point
		at the origin of a first-order linear system of differential equations.
		In order to understand the Stokes phenomenon they perturb the system,
		introducing a small parameter, which causes the splitting of the irregular
		singularity into two finite Fuchsian singularities. 
		Instead of considering first-order linear system, we study 
		a higher-order scalar equation. To make our discussion 
		simpler, we investigate a solvable (reducible) higher-order equation. 	
	   With the present article we begin a research on  the nature of the
		newly introduced perturbation. We start with the most simple example, namely,
			with the zero level perturbation of a certain third-order initial equation.
		The main goal of this paper is to show  that there exists a connection, 
		by a limit $\sqrt{\varepsilon} \rightarrow 0$,
			between  the analytic invariants at the origin (the Stokes matrices)	
		 of the initial equation and the analytic invariants
			around the finite resonant Fuchsian singularities (the monodromy matrices)
			of the perturbed equation. In general the monodromy matrices of the perturbed equation
			decompose into convergent and divergent parts (see \thref{monodromy}). 
			Consider for every fixed $\nu\in\RR$
			the particular sequence $\sqrt{\varepsilon}=\sqrt{\varepsilon_n}$  for
			$1/\sqrt{\varepsilon_n}=\nu+2 n$. It turns out that along this sequence (which in fact
			defines the so called logarithmic resonant cases) the divergent pats of the monodromy
			matrices stay constant. Thus along this very particular sequence of values of $\varepsilon$
			the monodromy matrices convergent. As a consequence of the decomposition theorem
			of Lambert and Rousseau (Proposition 4.31 in \cite{CL-CR}), a theorem of Klime\v s
			(Theorem 32. in \cite{MK3}) and our decomposition theorem (\thref{monodromy}) we 
			point which part of the monodromy matrices
			are the so called unfolded Stokes matrices. This result  allows us to connect by a limit 
			$\sqrt{\varepsilon} \rightarrow 0$ the pointed part of the monodromy matrices and
			the Stokes matrices of the initial equation (see \thref{main}).

		  		Problems of this kind have been considered in the already  classical works of 
		  		Glutsyuk \cite{AG}, Ramis \cite{R}, Zhang  \cite{Z}, Duval \cite{D}, Sch\"afke \cite{S} 
					as well as in recent ones: Glutsyuk \cite{AG1}, Lambert and Rousseau \cite{HLR, CL-CR, CL-CR1, CL-CR2},
					Slavyanov and Lay \cite{S-L}, Klime\v s \cite{MK, MK1, MK2, MK3}, Remy \cite{Re}.
						 Our work is closer to the works of Glutsyuk \cite{AG, AG1}
						and  Lambert and Rousseau \cite{CL-CR, CL-CR1} where the authors introduce a
						small parameter that splits the irregular point at the origin of a linear system
						(not a scalar equation) into two finite Fuchsian singularities. Then they study
						the confluence of the connection matrices and the monodromy matrices of the perturbed system to 
						the Stokes matrices of the original system. 
						The work of Glutsyuk treats the confluence on sectors in the parameter 
						$\varepsilon$-space, on which the regular singularities are non-resonant.
						He shows, that, generically,  the limit $\varepsilon \rightarrow 0$ of no
						product of monodromy matrices does give the Stokes matrices (Theorem 4.6 in \cite{AG1}).
						 The approach of Lambert, Rousseau and Hurtubise \cite{HLR, CL-CR, CL-CR1} of
						``mixed bases'' allows to treat the resonant values of $\varepsilon$ and
						may be used in studying of a scalar equation too. However, the calculation
						of the monodromy matrices of a linear system (not a scalar equation),
						relevant to the same fundamental matrix solution with the resonant values of $\varepsilon$
				    requires more complicated computations. We also note that 
						the approach in \cite{CL-CR} allows to make limits of monodromy along
						any particular sequence $\varepsilon= \varepsilon_n$, with
						$1/\varepsilon_n=1/\varepsilon_0 + 2 n$, such as in Duval \cite{D}.
						It should be of use for studying in case with $\varepsilon\in\CC$

     Our approach to the above results is different. We fully exploit the resonance and the appearance
		of logarithmic terms in the mixed basis of solutions of the perturbed equation in order to calculate
		the  monodromy matrices.
	This approach is assisted by the representation of this basis
		 in terms of iterated integrals depending only on
			the solutions of the equations  $L_j(\varepsilon)\,u=0$. 
		As a result  we can easy determine 
			the coefficients in front of the present logarithmic terms of the solutions
			of the perturbed equation just as the residues at the singular points 
			of the functions under integration.  Moreover, our approach allows us to point
			which part of the monodromy matrices  take the part of the so called unfolded Stokes matrices. 
			We believe that the developed technique could be useful
			in case of a general higher-order solvable scalar equation. It could be also
			applied in case of a linear system, whose fundamental matrix solutions
			are represented in the form of appropriate integrals.

  The motivation for the study of the exactly this initial equation comes from the
	investigation of the integrability of the Painlev\'e equations. In particular,
	when $\nu=1/2$ our  initial  equation 
 appears as the second normal variational equation
	of the Hamiltonian system $\mathcal{H}_{IV}(y, p, t, a, b)$ corresponding to
	the fourth Painlev\'e equation
	  \ben
		 \ddot{y}=\frac{1}{2 y}\,(\dot{y})^2 + \frac{3}{2}\,y^3 + 4 t\,y^2
		  + 2 (t^2 - a)\,y + \frac{b}{y}
		\een
 along a particular solution $y=p=0,\,  a=1,\,b=0$ \cite{TS}. 
In the same paper \cite{TS} we have proved that
the connected component $G^0$ of the unit element of the differential Galois group
of the initial equation  with $\nu=1/2$ is not Abelian using Stokes matrices.
Such solvable differential equations whose differential Galois group is a solvable algebraic
group can be found in our previous work on the
Painlev\'e V (again 1 irregular point at origin of Poincar\'e rank 1) \cite{TS1} and on the
Painlev\'e VI equation (Fuchsian differential equations) \cite{HS, TS2}.

This article is organized as follows. In the next section we build global
fundamental ``mixed'' matrices 
of the initial and the perturbed equations, with respect to which
  we are going to determine the corresponding analytic invariants. In section 3
we explicitly compute the formal monodromy and the Stokes matrices of the
initial equation. In section 4 we explicitly  calculate
the monodromy matrices around the finite Fuchsian singularities of the perturbed
equation during a resonance. In section 5  we
establish the main results of this paper.

   
	\section { Global solutions }
	
	In this section we will introduce the global fundamental matrices of
	the initial and the perturbed equations, with respect to which we
	will compute in the next two sections the corresponding analytic
	invariants. As we have announced in the introduction 
	we are going to use fundamental matrices different from the usual form.
	
	Considered in this paper equations are  very particular cases of a more
		general third-order solvable ODE 
		       \be\label{e}
	      L\,y=0\,,\quad
				L=L_3 \circ L_2 \circ L_1\,,						
						\ee	
	where $L_j,\,1 \leq j \leq 3$ are first-order  differential operators
	of the form
		 \be\label{form-g}
		     L_j=\partial + a_j(x) \qquad \textrm{with}\qquad
				a_j(x)\in\CC(x)\,,\quad  1 \leq j \leq 3 \,.	
		 \ee
			
  To introduce the fundamental matrices we at first reduce a scalar solvable equation
	\eqref{e}-\eqref{form-g}   to a 
	certain special linear system. We will call such a system, a system
	associated with the given scalar equation.
 Denote $\tilde{L}_j=L_3  \circ \cdots \circ L_j,\, 1 \leq j \leq 3$,
as $\tilde{L}_1=L$, and by $z_j(x),\, 1 \leq j \leq 3$ a solution of the equation
$\tilde{L}_j\,u=0$. The function $z_1(x)=y(x)$ is a solution of the equation $L\,y=0$. Then we have

 \bth{t1}
 The function $y(x)$ is a solution of the third-order equation \eqref{e}-\eqref{form-g}
if and only if the vector
$Y(x)=(y(x), L_1\, y(x), L_2 \circ L_1\, y(x))^{\tau}$ solves the system
  \be\label{system}
	  Y'(x)=A(x)\,Y(x)
	\ee
	with
	\be\label{A1}
	  A(x)=\left(\begin{array}{ccc}
		 -a_1(x)   &1        &0\\
		   0       &-a_2(x)  &1 \\
			 0       &0        &-a_3(x)
			          \end{array}
					\right)\,.			
	\ee
\ethe

 \proof
The proof is straightforward after the observation that $z_2(x)=L_1 y(x)$
and $z_3(x)=L_2 \circ L_1\, y(x)$.\qed

\bde{d}
 A $3 \times 3$ matrix $\Phi(x)$ is called a fundamental matrix of the scalar
third-order equation \eqref{e}-\eqref{form-g} if it is a fundamental matrix of a system, associated
with the same scalar equation.
\ede

 It turns out that both initial and perturbed equations  have a global
fundamental matrix, whose elements are expressed in terms of iterated integrals
depending only on  the solutions
of the equations $L_j\,u=0,\,1 \leq j \leq 3$.
Denote by $\Phi(x, \cdot)$, where the second argument is either $0$ or $\varepsilon$,
the fundamental matrix of the initial or the perturbed equation respectively.

 \bth{solution}
	  Both initial and perturbed equations admit a global fundamental matrix
	 $\Phi(x, \cdot)$ of the form
			\be\label{fms}
			 \Phi(x, \cdot)=\left(\begin{array}{ccc}
			   \Phi_1(x, \cdot)    &\Phi_{12}(x, \cdot)  &\Phi_{13}(x, \cdot)\\
				   0                   &\Phi_2(x, \cdot)     &\Phi_{23}(x, \cdot)\\
					 0                   &0                      &\Phi_3(x, \cdot)
					                        \end{array}
														\right)\,,			
			\ee
			where the diagonal elements $\Phi_j(x, \cdot),\,j=1, 2, 3$ are the solutions
			of the equations $L_j(\cdot) u=0,$ with
			$L_j(\varepsilon)$ and $L_j(0)=L_j$ given by \eqref{pe} and \eqref{npe} respectively.
			The off-diagonal elements are defined as follows,
			 \be\label{iint}
			       \Phi_{12}(x, \cdot) &=&
				\Phi_1(x, \cdot) \int_{\Gamma_1(x, \cdot)} \frac{\Phi_2(t_1, \cdot)}{\Phi_1(t_1, \cdot)}\,
						                   d t_1\,,\\[0.4ex]
					\Phi_{23}(x, \cdot)   &=&
					\Phi_2(x, \cdot) 
					    \int_{\Gamma_2(x, \cdot)} \frac{\Phi_3(t_2, \cdot)}{\Phi_2(t_2, \cdot)}\,
						                   d t_2\,,\nonumber\\[0.9ex]
				\Phi_{13}(x, \cdot)     &=&
				\Phi_1(x, \cdot) 
				   \int_{\Gamma_2(x, \cdot)} \frac{\Phi_2(t_1, \cdot)}{\Phi_1(t_1, \cdot)}
					    \left(\int_{\Gamma_2(t_1, \cdot)} \frac{\Phi_3(t_2, \cdot)}{\Phi_2(t_2, \cdot)}\,
							        d t_2\right) d t_1\,.\nonumber
         \ee
				 The paths of integration $\Gamma_j(x, \varepsilon)$ 
				and $\Gamma_j(x, 0)$  are taken from the same base point
				$x$ in such a way that $\Gamma_j(x, \varepsilon) \rightarrow \Gamma_j(x, 0)$
				as $\varepsilon \rightarrow 0\in\RR$, and the matrices $\Phi(x, \cdot)$
				are  fundamental matrix solutions of the initial and the perturbed equations
				respectively.
	\ethe

  \proof
	To prove the statement, we have to check that the so defined matrix $\Phi(x, \cdot)$ is such that
	$\Phi'(x, \cdot)=A(x)\,\Phi(x, \cdot)$, and $\det\,\Phi(x, \cdot) \neq 0$ outside of the singular points of the
	system \eqref{system}-\eqref{A1}.
	
	The first condition is checked directly. Note that
	$\det\,\Phi(x, \cdot)=\Phi_1(x, \cdot)\,\Phi_2(x, \cdot)\,\Phi_3(x, \cdot)$. 
	But for every function $\Phi_j(x, \cdot)$ 
	we have that $\Phi_j(x, \cdot)\neq 0$ outside of the singular points of the equation $L_j\,u=0$
	respectively, since it is a fundamental  solution of this equation. Now the second condition
	follows directly from the observations that the singular points of the system
	\eqref{system}-\eqref{A1} and of the initial and perturbed equations
	coincide.\qed

  As a direct corollary of \thref{solution} we obtain the following  proposition

  \bpr{p2}
	  Both initial and perturbed equations  possess
		a global fundamental set of solutions of the form
		 \ben
		     & &
		   		  \Phi_1(x, \cdot)\,,\,
			\Phi_1(x, \cdot) \int_{\Gamma_1(x, \cdot)} \frac{\Phi_2(t_1, \cdot)}{\Phi_1(t_1, \cdot)} d t_1\,,\\[0.4ex]
			   & &
			\Phi_1(x, \cdot) \int_{\Gamma_2(x, \cdot)} \frac{\Phi_2(t_1, \cdot)}{\Phi_1(t_1, \cdot)} 
			   \left(\int_{\Gamma_2(t_1, \cdot)}  \frac{\Phi_3(t_2, \cdot)}{\Phi_2(t_2, \cdot)}  d t_2\right) d t_1\,. 
			\een
	\epr


    \section{ The analytic invariants of the initial equation }

  In this section we will introduce and compute by hand the
 formal monodromy and the Stokes matrices at the origin of the initial equation.
	In this paper we are going to use the summability theory  (applied
	to ordinary differential equations) to calculate the Stokes matrices.

	All singular directions and sectors are defined on the Riemann surface of the natural logarithm.
 Consider the initial equation together with the formal fundamental matrix at the origin
(in the form of the theorem of Hukuhara-Turrittin \cite{W}), given by the next proposition
  \bpr{formal}
	  \,\,\,The initial equation admits a unique formal \,\,fundamental matrix
		$\hat{\Phi}(x, 0)$ at the origin in the form 
	    	\be\label{formals}
	  \hat{\Phi}(x, 0)=\hat{H}(x) \,x^{\Lambda}\,
		\exp\left(-\frac{Q}{x}\right)\,,
	\ee	
	where
		  \be\label{l}
	  Q=\diag(1,  2, 0)\,,\qquad
		\Lambda=\diag(0, \nu-2, \nu-4)
	\ee
	  and
	\be\label{H-f}
	  \hat{H}(x)=\left(\begin{array}{ccc}
		  1    &x^2\,\hat{\varphi}(x)      &\frac{x^4\,\hat{\psi}(x)}{2} \\[0.3ex]
			0    &1                     &-\frac{x^2}{2}\\[0.3ex]
			0    &0                     &1
			               \end{array}
							\right)\,.			
				\ee
		The elements $\hat{\varphi}(x)$ and $\hat{\psi}(x)$ of the matrix
				$\hat{H}(x)$ are defined as follows,
				\begin{enumerate}
				  \item\,
					 If $\nu$ is a non-positive integer, then 
					   \be\label{fin}
	      & &
	  \hat{\psi}(x)=1 + \nu\,x + \nu\,(\nu+1)\,x^2 + \nu\,(\nu+1)\,(\nu+2)\,x^3 + \cdots
		+ (-1)^{\nu}\,(-\nu)!\,x^{-\nu}\,,\\[0.2ex]
        & &
		\hat{\varphi}(x)=1 - \nu\,x + \nu\,(\nu+1)\,x^2 - \nu\,(\nu+1)\,(\nu+2)\,x^3 + \cdots
		+ (-\nu)!\,x^{-\nu}\nonumber		
	\ee
	 are analytic at the origin functions.
	   \item\,
		Otherwise, $\hat{\psi}(x)$ and $\hat{\varphi}(x)$ are represented in terms of divergent
		series
	      \be\label{for}
	  \hat{\psi}(x)  &=&
		1 + \nu\,x + \nu\,(\nu+1)\,x^2 + \nu\,(\nu+1)\,(\nu+2)\,x^3 + \cdots\,,\\[0.2ex]
		\hat{\varphi}(x) &=&
		1 - \nu\,x + \nu\,(\nu+1)\,x^2 - \nu\,(\nu+1)\,(\nu+2)\,x^3 + \cdots\nonumber
	\ee	 
				\end{enumerate}
	\epr
	
	  \proof
  	Choosing $\hat{\Phi}_1(x, 0)=e^{-1/x}\,,\,\hat{\Phi}_2(x, 0)=x^{\nu-2}\,e^{-2/x}\,,
		\,\,\hat{\Phi}_3(x, 0)=x^{\nu-4}$
	we obtain $\hat{\Phi}_{23}(x, 0)=-x^{\nu-2}/2$ where 
	 $\Gamma_2(x, 0)$ in \eqref{iint} is a path from $0-$ to $x$, approaching $0$ in the
	direction $\RR_-$. Next, looking for
	$\hat{\Phi}_{12}(x, 0)$ and $\hat{\Phi}_{13}(x, 0)$ in the form
	$\hat{\Phi}_{12}(x, 0)=x^{\nu}\,e^{-2/x}\,\hat{\varphi}(x)$ and
	$\hat{\Phi}_{13}(x, 0)=x^{\nu}\,\hat{\psi}(x)/2$ we find
  that $\hat{\varphi}(x)$ and $\hat{\psi}(x)$ satisfy the following first-order equations
	 \be\label{nhe}
	      	   x^2\,\hat{\psi}' + (\nu\,x - 1)\,\hat{\psi}=-1\,,\qquad
		 x^2\,\hat{\varphi}' + (\nu\,x + 1)\,\hat{\varphi}=1\,.
	 \ee
	For almost all values of the parameter $\nu\in\CC$ equations \eqref{nhe} admit 
	particular solutions in terms of divergent series. Only in the  
	exceptional cases when $\nu=0, -1, -2, -3, \ldots$ the functions
	 \ben
	      & &
	  \hat{\psi}(x)=1 + \nu\,x + \nu\,(\nu+1)\,x^2 + \nu\,(\nu+1)\,(\nu+2)\,x^3 + \cdots
		+ (-1)^{\nu}\,(-\nu)!\,x^{-\nu}\,,\\[0.2ex]
        & &
		\hat{\varphi}(x)=1 - \nu\,x + \nu\,(\nu+1)\,x^2 - \nu\,(\nu+1)\,(\nu+2)\,x^3 + \cdots
		+ (-\nu)!\,x^{-\nu}		
	\een
	 are  particular solutions of equations \eqref{nhe}, analytic at the origin.

	Let $\nu \neq 0, -1, -2, -3, \ldots$. Then the divergent power series
	 \ben
	  \hat{\psi}(x)  &=&
		1 + \nu\,x + \nu\,(\nu+1)\,x^2 + \nu\,(\nu+1)\,(\nu+2)\,x^3 + \cdots\,,\\[0.2ex]
		\hat{\varphi}(x) &=&
		1 - \nu\,x + \nu\,(\nu+1)\,x^2 - \nu\,(\nu+1)\,(\nu+2)\,x^3 + \cdots
	\een
	are particular solutions of the equations \eqref{nhe}.
	Now it is easy to write down the formal fundamental matrix $\hat{\Phi}(x, 0)$
	in the wished form.
	
	This completes the proof.\qed 
	
	 Now we can make the formal monodromy matrix $\hat{M}$ explicit.
		 \bpr{fmon}
  The formal monodromy matrix $\hat{M}$ relative to the formal solution \eqref{formals}
	 is defined by
	 		  $$\,
			  \hat{\Phi}(x.e^{2 \pi\,i}, 0)=\hat{\Phi}(x, 0).\hat{M}\,,
			\,$$
		where
	  \be\label{fm}
		  \hat{M}=e^{2 \pi\,i\,\Lambda}=\left(\begin{array}{ccc}
			  1    &0      &0\\
				0    &e^{2 \pi\,i\,\nu}   &0\\
				0    &0      &e^{2 \pi\,i\,\nu}
				            \end{array}
						\right)\,.				
		\ee
		\epr
 The formal monodromy $\hat{M}$ is a formal analytic invariants of the
	initial equation.
 
\bde{sector}
	   1.\,A sector is to be a set of the form
		 $$\,
		  S=S(\theta, \alpha, \rho)=\{\,x=r\,e^{i \delta}\,|\,
			0 < r < \rho,\,\theta-\alpha/2 < \delta < \theta + \alpha/2\,\}\,,
		\,$$
		where $\theta$ is an arbitrary real number (the bisector of $S$),
		$\alpha$ is a positive real (the opening of $S$) and $\rho$ is either 
		a positive number or $+ \infty$ (the radius of $S$).\\
		2.\, A closed sector is a set of the form
		 $$\,
		  \bar{S}=\bar{S}(\theta, \alpha, \rho)=\{\,x=r\,e^{i \delta}\,|\,
			0 < r \leq \rho,\,\theta-\alpha/2 \leq \delta \leq \theta+\alpha/2\,\}\,,
		\,$$
		with $\theta$ and $\alpha$ as before, but where $\rho$ is a positive real number
		(never equal to $+ \infty$).
	\ede

  From \prref{formal} it follows that the initial equation has singular directions,
	relative only to the divergent series $\hat{\psi}(x)$ and $\hat{\varphi}(x)$.
	Then the restriction of the definition of the singular directions \cite{M} in general
	gives us the following 
   \bde{ds}
	  For  the divergent series $\hat{\psi}(x)$ we define the 
	 singular direction (anti-Stokes direction) $\theta_1,\,0 \leq \theta_1 < 2 \pi$ 
	of the initial equation  as the bisector of the
  sector where $Re \left((-q_1+q_3)/x\right)=Re (-1/x) < 0$. Then $\theta_1=0$.
  In the same manner, the singular direction $\theta_2,\,0 \leq \theta_2 < 2 \pi$
	relative to the divergent series $\hat{\varphi}(x)$ is the bisector of the 
	sector where $Re ((-q_1+q_2)/x)=Re (1/x) < 0$. Then $\theta_2=\pi$.
   \ede

   Our next step is to make explicit an actual fundamental matrix $\Phi(x, 0)$
	at the origin, associate with the above formal fundamental matrix $\hat{\Phi}(x, 0)$
	via the next fundamental theorem
	 \bth{actual}(Hukuhara-Turrittin-Martinet-Ramis \cite{JM-JR, M})
 In the formal fundamental matrix $\hat{\Phi}(x, 0)$ \eqref{formals} of 
the 	initial  equation  defined by Hukuhara-Turrittin
 the entries of the matrix $\hat{H}(x)$ are 1-summable in every non-singular
direction $\theta$. If we denote by $H(x)$ the 1-sum of $\hat{H}(x)$ along
$\theta$ obtained from $\hat{H}(x)$ by a Borel - Laplace transform, then
$\Phi(x, 0)=H(x)\,x^{\Lambda}\,\exp\left(-Q/x\right)$
is an actual fundamental matrix of the 	initial equation. 
\ethe		
		
Let us briefly recall some definitions and facts needed to to build
1-sum (Borel sum) of the matrix $\hat{H}(x)$, following
the works of Balser \cite{WB}, Ramis \cite{R1, R2}, van der Put and Singer \cite{dPS}.

	\bde{gevrey}
	 A formal power series $\hat{f}(x)=\sum_{n=0}^{\infty} f_n\,x^n$ is said to be
	 of Gevrey order 1 if there exist two positive constants $C, A > 0$ such that
	  $$\,
		  |f_n| \leq C\,A^n\,n!\quad
			\textrm{for every}\quad
			 n\in\NN\,.
		\,$$
	\ede
	We denote by $\CC[[x]]_1$ the set of all power series of Gevrey order 1.
	
	\bde{borel}
	 The formal Borel transform $\hat{\B}_1$ of order 1 of a formal power series
	 $\hat{f}(x)=\sum_{n=0}^{\infty} f_n\,x^n$ is called the formal series
	  $$\,
		 (\hat{\B}_1\,\hat{f})(\zeta)=\sum_{n=0}^{\infty} \frac{f_n}{n!} \zeta^n\,.
		\,$$
	\ede
	If $\hat{f}\in\CC[[x]]_1$ then its formal Borel transform $(\hat{\B}_1\,\hat{f})$
	of order 1 converges in a neighborhood of the origin $\zeta=0$ with a sum $f(\zeta)$.
	
	 The inverse operator of the Borel transform is the Laplace transform.
	
	\bde{laplace}
	 Let $f(\zeta)$ be analytic and of exponential size at most 1 at $\infty$,
	i.e. $|f(\zeta)| \leq A\,\exp(B |\zeta|),\,\zeta\in \theta$ along any direction 
	$\theta$ from $0$ to $+\infty\,e^{i \theta}$. Then the integral
	 $$\,
	  (\L_\theta f)(x)=\int_0^{+\infty e^{i \theta}} f(\zeta)\,\exp\left(-\frac{\zeta}{x}\right)\,
		d\left(\frac{\zeta}{x}\right)
	\,$$
	is said to be the Laplace complex transform $\L_{\theta}$ of order 1 in the direction $\theta$ of $f$.
	\ede
	
	\bde{summable}
	  The formal power series $\hat{f}(x)=\sum_{n=0}^{\infty} f_n\,x^n$ is 1-summable
		(or Borel summable) in the direction $\theta$ if there exist an open sector $V$ 
		bisected by $\theta$ whose 
		opening is $> \pi$ and a holomorphic function $f(x)$ on $V$ such that for every
		non-negative integer $N$,
		  $$\,
			 \left| f(x) - \sum_{n=0}^{N-1} f_n\,x^n \right| \leq C_{V_1}\,A^N_{V_1}\,N!\,|x|^N
			\,$$
			on every closed subsector $\bar{V}_1$ of $V$ with constants $C_{V_1}, A_{V_1} > 0$
			depending only on $V_1$. The function $f(x)$ is called the 1-sum (or Borel sum) of $\hat{f}(x)$
			in the direction $\theta$.
	\ede
	If $\hat{f}(x)$ is 1-summable in all but a finite number of directions, we will say that
	it is 1-summable. 
	
	One useful criterion for a Gevrey series of order 1 to be 1-summable is given in terms
	of Borel and Laplace transforms.
	
	 \bpr{crit}$($\cite{dPS}$)$
	  Let $\hat{f}\in\CC[[x]]_1$ and let $\theta$ be a direction. The the following are
		equivalent:\\
		 1.\,$\hat{f}$ is 1-summable in the direction $\theta$.\\
		 2.\,The convergent power series $(\hat{\B}_1 \hat{f})(\zeta)$ has an analytic continuation
		$h$ in a full sector $\{\,\zeta\in\CC\,|\, 0 < |\zeta| < \infty,\,
		|\arg(\zeta) - \theta | < \epsilon\,\}$. In addition, this analytic continuation has
		exponential growth of order $\leq 1$ at $\infty$ on this sector, i.e.
		$|h(\zeta)| \leq A\,\exp(B |\zeta|)$. In this case $f=\L_{\theta}(h)$ is its 1-sum.
		  \epr

	Applying the above theory to the divergent power series 
 $\hat{\varphi}(x)$ and $\hat{\psi}(x)$ we obtain their 1-sums (Borel sums).
	
	\ble{l1}
Let $\nu \neq 0, -1, -2, -3, \ldots$. 
Then for any directions $ \theta_1 \neq 0$ and
$\theta_2 \neq \pi$ from 0 to $+\infty\,e^{i \theta_k},\,k=1, 2$  
the functions
 \be\label{sums}
       & &
  \psi_{\theta_1}(x)=x^{-1}\,\int_0^{+\infty e^{i \theta_1}} (1 - \zeta)^{-\nu}\,
	                 e^{-\frac{\zeta}{x}}\, d \zeta\,,\\[0.5ex]
			 & &						
	\varphi_{\theta_2}(x)=x^{-1}\,\int_0^{+\infty e^{i \theta_2}} (1 + \xi)^{-\nu}\,
									e^{-\frac{\xi}{x}}\,d \xi\nonumber
\ee
define the 1-sum (Borel sum) of the divergent series $\hat{\psi}(x)$ and 
$\hat{\varphi}(x)$, respectively,   in such directions.
\ele  

  \proof 
Let us represent the divergent series $\hat{\psi}(x)$ and $\hat{\varphi}(x)$ as
 $$\,
  \hat{\psi}(x)=\sum_{n=0}^{\infty} (\nu)^{(n)}\,x^n\,,\qquad
	\hat{\varphi}(x)=\sum_{n=0}^{\infty} (-1)^n\,(\nu)^{(n)}\,x^n\,,
\,$$
where  $(\nu)^{(n)}$ for $n=0, 1, 2, \ldots$ is the rising factorial
 $$\,
   (\nu)^{(n)}=\nu\,(\nu+1)\,(\nu+2) \ldots (\nu+n-1)\,,\quad
	(\nu)^{(0)}=1\,.
\,$$
Let $|\nu| \leq 1$. Then
$$\,
   |(\nu)^{(n)}| \leq |\nu|\,(|\nu|+1)\,(|\nu|+2) \ldots (|\nu|+n-1) \leq n!\,.
\,$$
Let $|\nu| > 1$. Then using the fact $1 < |\nu| + 1$, we have the following
estimates
  \ben
	  & &
	  |\nu| < |\nu|+1\,,\quad |\nu|+1 < |\nu|+1 + 1 < 2(|\nu+1)\,,\quad
		|\nu|+2 < 3 (|\nu|+1)\,, \ldots\\[0.2ex]
		& &
		|\nu|+n-1 < n (|\nu|+1)\,.
	\een
	Then 
	$$\,
	  |(\nu)^{(n)}| < (|\nu|+1)^n\,n!.
	\,$$
	We have the same estimates for $(-1)^n\,(\nu)^{(n)}$. Therefore the series
	$\hat{\psi}(x)$ and $\hat{\varphi}(x)$ are of Gevrey order 1 with constants
	$C=A=1$ if $|\nu| \leq 1$ and $C=1,\,A=|\nu|+1$ if $|\nu| > 1$ (see \deref{gevrey}).
	
	As a result their formal Borel transforms (\deref{borel})
	 \ben 
	   & &
	  (\hat{\B}_1 \hat{\psi})(\zeta)=\sum_{n=0}^{\infty}
		 \frac{(\nu)^{(n)}}{n!} \zeta^n=(1-\zeta)^{-\nu}\,,\\[0.4ex]
		& &
		(\hat{\B}_1 \hat{\varphi})(\xi)=\sum_{n=0}^{\infty}
		 \frac{(-1)^n\,(\nu)^{(n)}}{n!} \xi^n=(1+\xi)^{-\nu}
	\een
	are analytic functions near the origin in the Borel planes.
	
	Then for any directions $\theta_1 \neq 0$ and $\theta_2 \neq  \pi$ 
	from 0 to $+\infty\,e^{i \theta_k},\,k=1,2$, the associate Laplace transforms
	(\deref{laplace})
	 \ben
	    & &
	 \psi_{\theta_1}(x)=
	\int_0^{+\infty e^{i \theta_1}} (1-\zeta)^{-\nu}\,
	\exp\left(-\frac{\zeta}{x}\right) d\left(\frac{\zeta}{x}\right)\,,\\[0.7ex]
	   & &
	 \varphi_{\theta_2}(x)=
	\int_0^{+\infty e^{i \theta_2}} (1+\xi)^{-\nu}\,
	 \exp\left(-\frac{\xi}{x}\right) d\left(\frac{\xi}{x}\right)
	\een
	define the corresponding 1-sum of the series $\hat{\psi}(x)$ and $\hat{\varphi}(x)$
	respectively in such directions (\prref{crit}).
	
	This completes the proof.\qed

   \bre{borel}
	    Moving the direction $\theta_1$ (resp. $\theta_2$) continuously 
	the corresponding Borel sums $\psi_{\theta_1}(x)$ (resp. $\varphi_{\theta_2}(x)$)
	stick each other analytically and define an  function $\tilde{\psi}(x)$
	(resp. $\tilde{\varphi}(x)$) on a sector of opening 
	$3 \pi,\,-\pi/2 < \arg x < 5 \pi/2$
	(resp. $- 3 \pi/2 < \arg x < 3 \pi/2$). In these sectors the multivalued functions
	$\tilde{\psi}(x)$ and $\tilde{\varphi}(x)$ define the Borel sums of  the series $\hat{\psi}(x)$ 
	and $\hat{\varphi}(x)$, respectively. 
	In every non-singular direction $\theta$ the multivalued functions $\tilde{\psi}(x)$
	and $\tilde{\varphi}(x)$ have one value $\psi_{\theta}(x)$ and $\varphi_{\theta}(x)$,
	respectively. Near the singular direction $\theta=0$ the function $\tilde{\psi}(x)$
	has two different values: $\psi^+_0(x)=\psi_{0+\epsilon}(x)$ and
	$\psi^-_0(x)=\psi_{0-\epsilon}(x)$, where $\epsilon > 0$ is a small number. Similarly, near
	the singular direction $\theta=\pi$ the function $\tilde{\varphi}(x)$ has two different values:
	$\varphi^+_{\pi}(x)=\varphi_{\pi+\epsilon}(x)$ and $\varphi^-_{\pi}(x)=\varphi_{\pi-\epsilon}(x)$.
	\ere

   Replacing the elements $\hat{\varphi}(x)$ and $\hat{\psi}(x)$ of the matrix
	$\hat{H}(x)$ in \eqref{H-f} by their sums (classical when $\nu$ is a non-positive integer
	and Borel otherwise), we obtain an actual function $H(x)$. It together with the actual
	function $F(x)=x^{\Lambda}\,\exp(-Q/x)$ form an actual fundamental matrix of the
	initial equation at the origin.

   \bpr{actual}
	  \begin{enumerate}
		 \item\,
	 Assume that $\nu$ is a non-positive integer. Then the initial equation possesses
	an unique actual fundamental matrix $\Phi(x, 0)$ at the origin in the form
	 $$\,
	  \Phi(x, 0)=H(x)\,F(x)\,,
		\,$$
		where $H(x)=\hat{H}(x)$ is an analytic at the origin function, defined by \eqref{H-f}, whose
		elements are  the analytic functions \eqref{fin}. The matrix $F(x)$ is the branch of
		$x^{\Lambda}\,e^{-Q/x}$ for $\arg x$.
		  \item\,
		Assume that $\nu$ is not a non-positive integer. Then 
	for every non-singular direction $\theta$  the initial equation possesses 
	 an unique actual fundamental matrix $\Phi_{\theta}(x, 0)$ at the origin in the form
	  \be\label{ac}
		 \Phi_{\theta}(x, 0)= 
		H_{\theta}(x)\,F_{\theta}(x)\,,
		\ee
		where  $H_{\theta}(x)$ is the Borel sum of the matrix $\hat{H}(x)$ in this direction and $F_{\theta}(x)$ is the
		branch of $x^{\Lambda}\,e^{-Q/x}$ for $\arg x=\theta$. In particular, 
		$\Phi_{\theta+2 \pi}(x, 0)=\Phi_{\theta}(x, 0)\,\hat{M}$.
		
		For the singular direction $\theta=0$ the initial equation has two actual fundamental matrices
		 \be\label{ac-0}
		   	 \Phi^{\pm}_0(x, 0)=\Phi_{0\pm\epsilon}(x, 0) \,,
		\ee
	where the matrices $\Phi_{0\pm \epsilon}(x, 0)$ are given by \eqref{ac} for a small number $\epsilon > 0$.

		  For the singular direction $\theta=\pi$ the initial equation again has two actual fundamental matrices
		 \be\label{ac-1}
		   	 \Phi^{\pm}_{\pi}(x, 0)=\Phi_{\pi\pm\epsilon}(x, 0) \,,
		\ee
	where the matrices $\Phi_{\pi\pm \epsilon}(x, 0)$ are again given by \eqref{ac}. 
   \end{enumerate} 
Moreover, the matrix $\Phi(x, 0)$ defined by \eqref{ac}-\eqref{ac-0}-\eqref{ac-1}
	and the fundamental matrix $\Phi(x, 0)$ introduced by \eqref{fms}-\eqref{iint}
	define the same actual fundamental matrix at the origin of the initial equation. 
	\epr
	
	\proof
	From \thref{actual}, \prref{formal}, \leref{l1} and \reref{borel} it follows
	that the matrices $\Phi(x, 0)$ defined by \eqref{ac}-\eqref{ac-0}-\eqref{ac-1} 
	 are the unique actual fundamental matrices at the origin,
	associated with the pointed formal fundamental matrix $\hat{\Phi}(x, 0)$.

	Therefore, we have only to show that these matrices 
	 and the fundamental matrix $\Phi(x, 0)$, introduced
	by \eqref{fms} and \eqref{iint} define the same actual fundamental matrix solution
	at the origin of the initial equation.

Let us represent $\Phi_{12}(x, 0)$,
	from \eqref{iint}, in the following form
	 \ben
	  \Phi_{12}(x, 0)
		     &=&
		\Phi_1(x, 0)\,\int_0^x \frac{\Phi_2(t, 0)}{\Phi_1(t, 0)}\,d t=
		e^{-1/x}\,\int_0^x t^{\nu-2}\,e^{-1/t}\,d t=\\[0.4ex]
		     &=&
		e^{-2/x}\,\int_0^x t^{\nu-2}\,e^{-\frac{1}{t}+\frac{1}{x}}\, d t 
	\een
	where  $\Gamma_1(x, 0)$ in \eqref{iint} is a path from $0+$ to $x$, approaching $0$
	 in the direction $\RR_+$.
	Then, introducing a new variable $\xi$ via 
	$$\,
	-\frac{1}{t}+\frac{1}{x}=-\frac{\xi}{x}
	\,$$
	we obtain
	 \ben
	  \Phi_{12}(x, 0)=x^{\nu-1}\,e^{-2/x}\,
		\int_0^{+\infty} (1+\xi)^{-\nu}\,e^{-\frac{\xi}{x}}\,d \xi\,.
	\een
  In the same manner, we can represent  $\Phi_{13}(x, 0)$ from \eqref{iint} as
	\ben
	  \Phi_{13}(x, 0)=\Phi_1(x, 0)\int_{\Gamma_2(x, 0)} \frac{\Phi_{23}(t, 0)}{\Phi_1(t, 0)}\, d t=
		-\frac{1}{2}\,\int_0^x
		 t^{\nu-2}\,e^{\frac{1}{t}-\frac{1}{x}}\,d t\,,
	\een
	where the path of integration is the path $\Gamma_2(x, 0)$ from $0-$ to $x$,
	used in the definition of $\hat{\Phi}_{23}(x, 0)$. 
Again, by  introducing  a new variable $\zeta$ by
$$\,
\frac{1}{t}-\frac{1}{x}=-\frac{\zeta}{x}\,,
\,$$
 we obtain the function
  \ben
	  \Phi_{13}(x, 0)=\frac{x^{\nu-1}}{2}\,
		 \int_0^{-\infty} (1-\zeta)^{-\nu}\,e^{-\frac{\zeta}{x}}\,d \zeta\,. 
	\een
Analytic continuations of the so built $\Phi_{12}(x, 0)$ and
$\Phi_{13}(x, 0)$ on $x$-plane yield  analytic functions
 \ben
   (\Phi_{12}(x, 0))_{\theta}= 
		x^{\nu}\,e^{-2/x}\,\varphi_{\theta}(x)\,,\quad
	(\Phi_{13}(x, 0))_{\theta}=
\frac{x^{\nu}}{2}\,\psi_{\theta}(x)
\een
on every non-singular direction $\theta$.

Let $\epsilon > 0$ be a small number and let $0+\epsilon$ and $0-\epsilon$
be two non-singular directions near the singular direction $\theta=0$. Then
the function $\Phi_{13}(x, 0)$ has two different values near $\theta=0$
 \ben
(\Phi_{13}(x, 0))^{\pm}_0=(\Phi_{13}(x, 0))_{0\pm \epsilon}=
\frac{x^{\nu}}{2}\,\psi_{0\pm \epsilon}(x)\,,
\een
where $\psi_{0 \pm \epsilon}(x, 0)$ are the Borel sums of the series $\hat{\psi}(x)$, built
by \leref{l1} and extended by \reref{borel}.

Similarly, let $\pi+\epsilon$ and $\pi-\epsilon$ be two non-singular directions
near the singular direction $\theta=\pi$. Then near $\theta=\pi$ the function 
$\Phi_{12}(x, 0)$ has two different values
 \ben
   (\Phi_{12}(x, 0))^{\pm}_{\pi}=(\Phi_{12}(x, 0))_{\pi\pm \epsilon}=
	x^{\nu}\,e^{-2/x}\,\varphi_{\pi \pm \epsilon}(x)\,.
\een 
Here $\varphi_{\pi \pm \epsilon}(x)$ are the Borel sums of the series $\hat{\varphi}(x)$,
built and extended by \leref{l1} and \reref{borel}.

In the same manner near $\theta=2 \pi$ the function $\Phi_{13}(x, 0)$ has two
different values
 \ben
  (\Phi_{13}(x, 0))^{\pm}_{2 \pi}=\frac{x^{\nu}\,e^{2 \pi\,i\,\nu}}{2}\,
	\psi_{2\pi \pm \epsilon}(x)\,.
\een
Note that near $\theta=2 \pi$ the actual function $F(x)$ is transformed  into
the function $F(x)\,\hat{M}=x^{\Lambda}\,e^{-Q/x}\,\hat{M}$.

This proves \prref{actual}.\qed

In what follows we define and compute the Stokes matrices of the initial equation.

Let $\theta$ be a singular direction of the initial equation. Denote by  $\Phi^+_{\theta}(x, 0)$ 
	and $\Phi^-_{\theta}(x, 0)$ 
	the actual fundamental matrix of the 	initial equation, defined by \prref{actual}. 
 Then
	
	\bde{stokes}
	 With respect to the given actual fundamental matrices the
	Stokes matrix $St_{\theta}\in \GL_3(\CC)$ corresponding to the singular direction $\theta$
	is defined by 
	  $$\,
	St_{\theta}=(\Phi^+_{\theta}(x, 0))^{-1}\,\Phi^-_{\theta}(x, 0)\,.
	 \,$$
	\ede

This definition implies that  the Stokes matrix measures the difference between two
fundamental matrices when we turn around a singular direction in a positive sense,
which is keeping with the definition of the monodromy matrices $M_j(\varepsilon), j=L, R$
of the perturbed equation (see next section).

  \bth{Stokes}
		  With respect to the formal fundamental matrix $\hat{\Phi}(x, 0)$ at
				the origin given by \eqref{formals}-\eqref{l}-\eqref{H-f}
				and associated to it the actual fundamental matrix at the origin given by
				\eqref{ac}-\eqref{ac-0}-\eqref{ac-1}, the initial equation has
				 two singular directions $\theta_1=0$ and $\theta_2=\pi$. The corresponding Stokes
					matrices are given by
					  				\ben 
				St_0=\left(\begin{array}{ccc}
     1     &0            &-\frac{\pi\,i}{\Gamma(\nu)}\\
		 0     &1            &0\\
		 0     &0            &1
						\end{array}
			\right)\,,\qquad							
				St_{\pi}=\left(\begin{array}{ccc}
	  1     &-\frac{2 \pi\,i\,e^{-\pi\,i\,\nu}}{\Gamma(\nu)}    &0\\
		0     &1                               &0\\
		0     &0                               &1
		             \end{array}
						\right)\,.							
\een
		\ethe

   \proof
	
The direction $\theta_1=0$ is a singular direction only to
the element $\Phi_{13}(x, 0)$ of the fundamental matrix. Therefore to
build the Stokes matrix $St_0$, corresponding to the singular direction
$\theta_1=0$, we have only to compare the functions 
$(\Phi_{13}(x, 0))^-_0=x^{\nu}\,\psi_{0-\epsilon}(x)/2$
and $(\Phi_{13}(x, 0))^+_0=x^{\nu}\,\psi_{0+\epsilon}(x)/2$.
We have  
\ben
 (\Phi_{13}(x, 0))^-_0=(\Phi_{13}(x, 0))^+_0+
\frac{x^{\nu-1}}{2}\int_\gamma (1-\zeta)^{-\nu}\,e^{-\frac{\zeta}{x}}\,d\,\zeta\,,
\een
where $\gamma=(0-\epsilon) - (0+\epsilon)$.
In this case, without changing the integral,  we can deform
$\gamma$  in a Hankel type path $\gamma'$ 
going along the positive real axis  from infinity to 1,
encircles 1 in the positive direction  and back to infinity in the positive sense. 
Then since $\arg(1-\zeta)=-\pi$ when $Re (\zeta) > 1$ and $\zeta$ lies on the direction 
$0+\epsilon$ the integral becomes  
    \ben
		      & &
		  \frac{x^{\nu-1}  (e^{-\pi\,i \nu}-e^{\pi\,i\,\nu})}{2}
			\int_1^{+\infty} (\zeta-1)^{-\nu} e^{-\frac{\zeta}{x}}\,d \zeta=\\[0.4ex]
			   &=&
			\frac{x^{\nu-1}\,e^{-1/x}( e^{-\pi\,i\,\nu}-e^{\pi\,i\,\nu})}{2}
			\int_0^{+\infty} u^{-\nu}\,e^{-\frac{u}{x}}\,d\,u=\\[0.4ex]
			   &=& 
				\frac{e^{-1/x} (e^{-\pi\,i\,\nu}-e^{\pi\,i\,\nu})}{2}
				\int_0^{+\infty} \tau^{-\nu}\,e^{-\tau}\,d \tau=
				\frac{1}{2} ( e^{-\pi\,i\,\nu}- e^{\pi\,i\,\nu})\,\Gamma(1-\nu)\,e^{-1/x}=\\[0.4ex]
				 &=&
				-\frac{\pi\,i}{\Gamma(\nu)}\,e^{-1/x}\,,
		\een
			where we used that the Gamma function $\Gamma(1-\nu)$ is related to $\Gamma(\nu)$  by
			(see \cite{HB-AE}) 
						$$\,
			  \Gamma(\nu)\,\Gamma(1-\nu)=\frac{\pi}{\sin(\pi\,\nu)}\,.
			\,$$

In the same manner, comparing the functions $(\Phi_{12}(x, 0))^-_{\pi}$ and
$(\Phi_{12}(x, 0))^+_{\pi}$, we obtain     
		\ben
			(\Phi_{12}(x, 0))^-_{\pi}=(\Phi_{12}(x, 0))^+_{\pi} -
			\frac{2 \pi\,i\,e^{-\pi\,i\,\nu}}{\Gamma(\nu)}\,e^{-1/x}\,.
		\een

 Then the straightforward application of the  \deref{stokes} gives us the
Stokes matrices at the origin of the initial equation 
\ben
  St_{\pi}=\left(\begin{array}{ccc}
	  1     &-\frac{2 \pi\,i\,e^{-\pi\,i\,\nu}}{\Gamma(\nu)}    &0\\
		0     &1                               &0\\
		0     &0                               &1
		             \end{array}
						\right)\,,\,\,\,\,\,\,\,
 St_0=\left(\begin{array}{ccc}
     1     &0            &-\frac{\pi\,i}{\Gamma(\nu)}\\
		 0     &1            &0\\
		 0     &0            &1
						\end{array}
			\right)\,.							
\een
We note that the function $1/\Gamma(\nu)$ is an entire function with zeros
at $\nu=0, -1, -2, \ldots$. Then, by a theorem on the analytic dependence of 
the Stokes matrices on the parameter $\nu$ (see \cite{HLR}), it follows 
(according to expectation), that  when $\nu\in\ZZ_{\leq 0}$ we have
$St_0=St_{\pi}=I_3$.

This completes the proof.
\qed

  Ir order to use the results of Hurtubuse, Lambert and Rousseau \cite{HLR, CL-CR, CL-CR1}, we
  consider now the initial equation and its actual fundamental matrices at the
	origin in the ramified domain 
		$\{x\in\CC\,:\, - \alpha < \arg (x) <  \alpha\}$, where $0 < \alpha < \pi/2$. 
		We cover this domain by 
   two open sectors $\Omega_1$ and $\Omega_2$ 
	  \ben
		  	\Omega_1 &=& \Omega_1(\alpha, \rho)=\left\{
			x=r\,e^{i \delta}\,|\, 0 < r < \rho,\, -\alpha < \delta < \pi + \alpha\right\}\,,\\[0.15ex]
		  \Omega_2 &=& \Omega_2(\alpha, \rho)=\left\{
			x=r\,e^{i \delta}\,|\, 0 < r < \rho,\, -(\pi+\alpha) < \delta < \alpha\right\}\,.
			\een
	 Denote by $\Omega_R$ and $\Omega_L$ the 
		connected components of the intersection $\Omega_1 \cap \Omega_2$. The radius $\rho$
		of the sectors $\Omega_j$ is chosen in such a way that $x_R=\sqrt{\varepsilon}$
		belongs to $\Omega_R$, and $x_L=-\sqrt{\varepsilon}$ belongs to $\Omega_L$.
		The lower sector $\Omega_2$ contains only the Stokes ray $i\,\RR_-$, 
		and the upper sector $\Omega_1$ contains
		only the Stokes ray $i\,\RR_+$. Then from the sectorial normalization theorem of
		Sibuya \cite{Sib} it follows, that over the sector $\Omega_1$ the matrices
		$H^+_0(x)$ and $H^-_{\pi}(x)$ represent the same actual function, asymptotic at the 
		origin in $\Omega_1$ to the matrix $\hat{H}(x)$ in \eqref{H-f}. Denote by
		$H_1(x)=H^+_0(x)=H^-_{\pi}(x)$.
		Similarly, over the sector
		$\Omega_2$ the matrices $H^+_{\pi}(x)$ and $H^-_0(x)$ represent the same analytic function,
		asymptotic at the origin in $\Omega_2$ to the matrix $\hat{H}(x)$. Denote by
		$H_2(x)=H^+_{\pi}(x)=H^-_0(x)$. Then the matrix
		 $$\,
		  \Phi_j(x)=H_j(x)\,F(x)\,\quad
			j=1, 2
		\,$$
		with the corresponding branch of $F(x)$ is an actual fundamental 
		matrix at the origin of the initial equation over the sector $\Omega_j,\,j=1, 2$.
	  Note that we can observe the Stokes phenomenon on
		$\Omega_R$ and $\Omega_L$. Let us turn around the origin in the positive sense,
		starting from the sector $\Omega_1$. On $\Omega_L$ we define the Stokes matrix
		$St_L$ as
		  $$\,
			  St_L=(\Phi_2(x, 0))^{-1}\,\Phi_1(x, 0)=St_{\pi}\,,
			\,$$
		where $St_{\pi}$ is the Stokes matrix, corresponding to the singular direction
		$\theta=\pi$ and defined by \thref{Stokes}. On $\Omega_R$ we define the Stokes
		matrix $St_R$ as
		  $$\,
			  (\Phi^+_0(x, 0))^{-1}\,\Phi^-_{2 \pi}(x, 0)=
				  (\Phi_1(x, 0))^{-1}\,\Phi_2(x, 0)\,\hat{M}=St_R\,\hat{M}=St_0\,\hat{M}\,,
			\,$$
		where $St_0$ is the Stokes matrix, corresponding to the singular direction
		$\theta=0$ and defined by \thref{Stokes}.
		Then the actual monodromy matrix $M_0$ around $x=0$ with respect to a base point on the lower 
		sector $\Omega_2$ and the corresponding fundamental solution there is given by
	    $$\,
			M_0=St_{\pi}\,St_0\,\hat{M}\,.
			\,$$ 
			Then the monodromy around $x=\infty$ is given by
	  $M_{\infty}^{-1}=St_{\pi}\,St_o\,\hat{M}$.
		In the last section we will make a cut between the singular points $x_L$ and $x_R$
		of the perturbed equation. Then the sectors $\Omega_j$ of this section will become
		new sectors $\Omega_j(\varepsilon)$ such that $\Omega_j(\varepsilon)$ tend to $\Omega_j$ when
		$\varepsilon$ tends to zero. This  limit procedure implies a limit between the
		monodromy matrices $M_j(\varepsilon)$ and the so called unfolded Stokes matrices $St_j(\varepsilon)$
		of the perturbed equation. The latter depend analytically on $\varepsilon$ and tend to
		the Stokes matrices $St_j, j=L, R$ defined above (see Theorem 4.25 in \cite{CL-CR}).

We end this section fixing the behavior of the other singularity of the initial equation.

\bre{E1}
  Except for the case when $\nu=0$ the point $x=\infty$ is a regular singular point for 
	the initial equation. The characteristic exponents $\rho^{\infty}_i,\,i=1, 2, 3$
	at $x=\infty$ are
	$$\,
	  \rho^{\infty}_1=0,\,\,\,\rho^{\infty}_2=1-\nu,\,\,\, \rho^{\infty}_3=2-\nu\,.
	\,$$
	When $\nu=0$ the change $x=1/t$ transforms the initial equation
	into equation
	 $$\,
	  \dddot{y}(t) + 3\,\ddot{y}(t) + 2\,\dot{y}(t)=0\,,\,\,\,\,\,\,\,\cdot=\frac{d}{d t}
	\,$$
	for which the point $t=0$ (resp. the point $x=\infty$ for the original equation) is an
	ordinary point.
\ere


  \section{The analytic invariants of the perturbed equation }

      In this section we will introduce and find the monodromy matrices of the perturbed equation,
			connecting it with the initial equation partially. 
 To do this, we firstly make the 
	global	fundamental matrix $\Phi(x, \varepsilon)$  \eqref{fms} of the perturbed equation
	explicit. Let us defines the paths of integrations $\Gamma_j(x, \varepsilon)$ in \eqref{iint}.
	Once fixing the paths of integration $\Gamma_j(x, 0)$, we 
	immediately determine the paths $\Gamma_j(x, \varepsilon)$ as follows:
	the path $\Gamma_1(x, \varepsilon)$ (resp. $\Gamma_2(x, \varepsilon)$) is a positive
	(resp. negative ) real trajectory of the vector field $(x^2 - \varepsilon)\,\partial_x$
	from $x_R=\sqrt{\varepsilon}$ to $x$ (resp. from $x_L=-\sqrt{\varepsilon}$ to $x$).
	The path $\Gamma_2(t_1, \varepsilon)$, similar to the path $\Gamma_2(t_1, 0)$, 
	is a path from $-\sqrt{\varepsilon}$ to $t_1$ in the direction $\RR_-$.
	Then we have that $\Gamma_j(x, \varepsilon) \rightarrow \Gamma_j(x, 0)$
	when $\varepsilon \rightarrow 0$. This choice of the paths implies that the parameter
	of perturbation $\varepsilon$ is a small real positive number, i.e.
	$0 < \varepsilon < 1$. 	
	 
	Next,  the perturbed equation is invariant under transformation
		 \be\label{symm}
		    \sqrt{\varepsilon} \longrightarrow - \sqrt{\varepsilon}\,.
		\ee
		So, throughout this section, we assume that $1/\sqrt{\varepsilon} > 1$.
		 
 The next readily verified Lemma simplifies the elements $\Phi_{23}(x, \varepsilon)$
and $\Phi_{13}(x, \varepsilon)$ of the fundamental matrix $\Phi(x, \varepsilon)$.		

\ble{int2}
  Let $a, b\in\RR$ such that $a > 0$ and $b>1$. Then
	  $$\,
		  \int_{-a}^x \frac{(s+a)^{b-1}}{(s-a)^{b+1}}\, d s=
			-\frac{1}{2 a\,b}\left(\frac{x+a}{x-a}\right)^b\,,
		\,$$
	where the integral is taken in the direction $\RR_-$ from $-a$ to $x < -a$.
\ele

 \thref{solution} in a combination with \leref{int2} gives the explicit
form of the fundamental matrix of the perturbed equation.
	 
	\bth{pe-fms}
	  Assume that $1/\sqrt{\varepsilon} > 1$.
		Then the explicit form of the elements of the fundamental matrix $\Phi(x, \varepsilon)$  
	is given	as follows
		  	\be\label{lsol}\quad
		  \Phi_1(x, \varepsilon)=\left(\frac{x-\sqrt{\varepsilon}}{x+\sqrt{\varepsilon}}\right)^{\frac{1}{2 \sqrt{\varepsilon}}}\,,
			\quad& &
		\Phi_2(x, \varepsilon)=(x^2-\varepsilon)^{\frac{\nu-2}{2}}\,\left(
			                     \frac{x-\sqrt{\varepsilon}}{x+\sqrt{\varepsilon}}
													                    \right)^{\frac{1}{\sqrt{\varepsilon}}}\,,\\[0,7ex]
							& &																
			\Phi_3(x, \varepsilon)=(x^2-\varepsilon)^{\frac{\nu-4}{2}}\nonumber																				
		\ee
		for the diagonal elements, and
   \be\label{lsol1}
	     & &
	   \Phi_{12}(x, \varepsilon)=\Phi_1(x, \varepsilon)\int_{\Gamma_1(x, \varepsilon)}
		  \frac{(t-\sqrt{\varepsilon})^{\frac{1}{2 \sqrt{\varepsilon}}+\frac{\nu-2}{2}}}
			     {(t+\sqrt{\varepsilon})^{\frac{1}{2 \sqrt{\varepsilon}}-\frac{\nu-2}{2}}}\,d t\,,\\[1.2ex]
				& &	 
				\Phi_{23}(x, \varepsilon)=-\frac{\Phi_2(x, \varepsilon)}{2}
						\left(\frac{x+\sqrt{\varepsilon}}{x-\sqrt{\varepsilon}}\right)^{\frac{1}{\sqrt{\varepsilon}}}=
						-\frac{1}{2}(x^2-\varepsilon)^{\frac{\nu-2}{2}}\,,\nonumber\\[1.2ex]
				& &	
				\Phi_{13}(x, \varepsilon)=-\frac{\Phi_1(x, \varepsilon)}{2}
			   \int_{\Gamma_2(x, \varepsilon)} 
				\frac{(t+\sqrt{\varepsilon})^{\frac{1}{2 \sqrt{\varepsilon}}+\frac{\nu-2}{2}}}
				     {(t-\sqrt{\varepsilon})^{\frac{1}{2 \sqrt{\varepsilon}}-\frac{\nu-2}{2}}}\,d t\,,\nonumber
			\ee
	for the off-diagonal elements.
	The path of integration $\Gamma_1(x, \varepsilon)$ is a path
	from $x_R=\sqrt{\varepsilon}$ to $x$ in the direction $\RR_+$. The path of integration
	$\Gamma_2(x, \varepsilon)$ is a path from $x_L=-\sqrt{\varepsilon}$ to $x$ in the direction
	$\RR_-$.
	 \ethe

   Let us briefly introduce the monodromy matrices of the perturbed equation, following
	Iwasaki et al.\cite{IKSY}.
  The perturbed equation is a third-order Fuchsian equation with three regular points
	over $\CC\PP^1$: two of them $x_R=\sqrt{\varepsilon}$ and $x_L=-\sqrt{\varepsilon}$ are 
	finite singularities and the third one is the infinity point.
	Let us consider the perturbed equation  over
		$X=\CC\PP^1 - \{x_R, x_L, x=\infty\}$.
  Its fundamental matrix $\Phi(x, \varepsilon)$  is a
   multi-valued analytic function on the punctured Riemann sphere $X$.
		Its multivaluedness  is described by the monodromy matrices.
		 Let $\gamma_j\in X,\,j=R, L$ 
		be  simple closed loops starting and ending at point $x_0=0\in X$, defined by
		   \be\label{loops}
				 \gamma_R(t)=\sqrt{\varepsilon} + \sqrt{\varepsilon}\,e^{\pi\,i(2 t+1)}\,,\quad
				 \gamma_L(t)=-\sqrt{\varepsilon} + \sqrt{\varepsilon}\,e^{2 \pi\,i\,t}\,,\quad
				 0 \leq t \leq 1\,.
				\ee
		Let the matrix $\Phi_{\gamma_j}(x, \varepsilon)$ be the analytic continuation of
		the fundamental matrix $\Phi(x, \varepsilon)$ of the perturbed equation
		along the loop $\gamma_j$. The matrix 
		 $\Phi_{\gamma_j}(x, \varepsilon)$ 
		depends only on the homotopy class $[\gamma_j]$ of $\gamma_j$.
		Since the perturbed equation is a linear equation, the matrix $\Phi_{\gamma_j}(x, \varepsilon)$ 
		is also a fundamental matrix of the same  equation and there is
		 a unique  invertible constant matrix $M_{\gamma_j}(\varepsilon)\in \GL_3(\CC)$ such that
			\be\label{M}
			  \Phi_{\gamma_j}(x, \varepsilon)=\Phi(x, \varepsilon)\,M_{\gamma_j}(\varepsilon)\,.
			\ee
    
      \bde{mr}
		 The antihomomorphism mapping  
				 \ben
				  & &
				  \pi_1(X, x_0) \longrightarrow \GL(3, \CC)\,,\\
					& &
				     [\gamma_j] \longrightarrow M_{\gamma_j}(\varepsilon)\,,\\
					& &	
						M_{\gamma_R\,\gamma_L}(\varepsilon)=M_{\gamma_L}(\varepsilon)\,M_{\gamma_R}(\varepsilon)\,,\quad
						M_{\gamma^{-1}_j}(\varepsilon)=M^{-1}_{\gamma_j}(\varepsilon)
				 \een
			determines monodromy representation of the perturbed equation 
			associated with the given fundamental matrix \cite{IKSY}.
			\ede
			The product $(\gamma_L\,\gamma_R)^{-1}$ of the so chosen loops $\gamma_j$ is homotopic
			to a simple loop $\gamma_{\infty}$ around infinity starting and ending at point $x_0$.
			Therefore the loops  $\gamma_j\,,\,j=R, L$ generate $\pi_1(X, x_0)$.
			  \bde{monodromy}
			The images $M_j(\varepsilon)=M_{\gamma_j}(\varepsilon)$ of the generators $\gamma_j, j=R, L$
			of $\pi_1(X, x_0)$ are called monodromy matrices of the perturbed equation  \cite{IKSY}.
			  \ede
		 They satisfy the following relation :
			\ben
				  M_L(\varepsilon)\,M_R(\varepsilon)=M^{-1}_{\infty}(\varepsilon)\,,
				\een
				where the matrix $M_{\infty}(\varepsilon)\in\GL(3, \CC)$ is the image of the loop $\gamma_{\infty}$.

	Now we have to determine when the matrix $\Phi(x, \varepsilon)$
	 contains logarithmic therms near the singular points
	$x_j, \,j=R, L$. In fact, only the elements $\Phi_{12}(x, \varepsilon)$ and
	$\Phi_{13}(x, \varepsilon)$ can contain such terms. To solve this problem,
	we are going to apply the local theory of the scalar Fucshian equations.
	Let us briefly recall some aspects of this theory needed to our goal, following 
	Golubev \cite{Go} and Iwasaki et al.\cite{IKSY}.
	For simplicity we restrict ourselves   to the perturbed equation.

As a Fuchsian equation, the perturbed equation can be written down as, \cite{Go}
	 \be\label{ee}
	 y'''(x)   &+&
	\frac{Q_2(x)}{(x-x_R)(x-x_L)}\,y''(x) + 
	  \frac{Q_1(x)}{(x-x_R)^2 (x-x_L)^2}\,y'(x) +\\[0.4ex]
		     &+&
		\frac{Q_0(x)}{(x-x_R)^3 (x-x_L)^3}\,y(x)=0\,,\quad
	  '=\frac{d}{d x}\,,\nonumber
	\ee
	where $Q_i(x),\,i=0, 1, 2$ are  polynomials of degree $3-i$. Denote by
	$c_i(x)$ the coefficient $Q_i(x)/((x-x_R)^{3-i} (x-x_L)^{3-i}),\,i=0, 1, 2$
	of the perturbed equation. At every regular point
		of the perturbed equation one can consider the so called characteristic 
		(or the inditial) equation.

	 \bde{ch}(\cite{Go})\,1.\,
	  The third order algebraic equation
		 $$\,
		  \rho\,(\rho-1)\,(\rho-2) + b_2\,\rho\,(\rho-1) + b_1\,\rho + b_0=0\,,
		\,$$
		where $b_i=\lim_{x\rightarrow x_j} c_i(x) (x-x_j)^{3-i},\,j=R, L$,
		is called the characteristic (or the indicial)  equation of the perturbed equation  at
		the regular singular point $x_j,\,j=R, L$. Its roots $\rho^j_i,\,j=R, L,\,i=1, 2, 3$
		are called the characteristic exponents at the singularities $x_j, \,j=R, L$.\\
		2.\, The characteristic equation at the point $t=0$ of the equation, obtained from
		the perturbed equation after the transformation $x=1/t$, is called the characteristic
		equation at the point $x=\infty$. Its roots $\rho^{\infty}_i,\,i=1, 2, 3$ are called
		the characteristic exponents at the regular point $x=\infty$
		
	\ede
	
	 It turns out that the coefficients $a_j(x, \varepsilon)\in\CC(x),\,j=1, 2, 3$ \eqref{pco}
	of the perturbed equation in the representation \eqref{perturbed}-\eqref{form}
	are expressed in the terms of the characteristic exponents.
	  \bpr{coeff}
		  The coefficients $a_j(x, \varepsilon)$ from \eqref{pco} of the 
			perturbed equation are unique-ly determined only by the
			characteristic exponents $\rho^j_i,\,j=R, L,\,i=1, 2, 3$ at the finite singularities
			$x_R$ and $x_L$.
		\epr
		
		\proof
	
	The coefficients $a_j(x, \varepsilon)$ must have the form:
	 \be\label{eq1}
	  a_1(x, \varepsilon)=-\frac{\rho^R_1}{x-x_R} - \frac{\rho^L_1}{x-x_L}\,,\quad
		a_2(x, \varepsilon)=-\frac{\rho^R_2-1}{x-x_R} - \frac{\rho^L_2-1}{x-x_L}\,,\\[0.7ex]
		a_3(x, \varepsilon)=-\frac{\rho^R_3-2}{x-x_R} - \frac{\rho^L_3-2}{x-x_L}\,.\nonumber
	  \ee
	\qed
	
	 According to \prref{coeff}, equations \eqref{eq1} and \eqref{pe}, the fundamental matrix of the 
	perturbed equation, as well as the monodromy representation with respect to this matrix are uniquely
		determined only by the characteristic exponents $\rho^j_i,\,j=R, L,\,i=1, 2, 3$
		$$\,
		 \rho^R_1=\frac{1}{2 \sqrt{\varepsilon}}\,,
		  \rho^R_2=\,\frac{\nu}{2} + \frac{1}{\sqrt{\varepsilon}}\,,\,
			\rho^R_3=\frac{\nu}{2};\qquad
		 \rho^L_1=-\frac{1}{2 \sqrt{\varepsilon}}\,,
		 \rho^L_2=\frac{\nu}{2} - \frac{1}{\sqrt{\varepsilon}}\,,\,
		\rho^L_3=\frac{\nu}{2}
		\,$$
	 at  the finite singular points $x_R=\sqrt{\varepsilon}$ and $\,x_L=-\sqrt{\varepsilon}$. 
	After the change $x=1/t$ the perturbed equation becomes
	 \be\label{inf}
	   \dddot{y}(t) 
		   &+& 
			\left(\frac{2 \nu}{t} + 3\right)\,\ddot{y}(t) +
		\left(\frac{\nu(\nu-1)}{t^2} + \frac{4\,\nu}{t} + 2\right)\,\dot{y}(t)+\\[0.3ex]
		   &+&
		\left(\frac{\nu\,(\nu-1)}{t^2} +
		\frac{2\,\nu}{t}\right)\,y(t)=0\,.\nonumber
	\ee
	According to \deref{ch}(2) the characteristic equation at $x=\infty$ of the perturbed equation
 is just the characteristic equation at $t=0$ of the equation \eqref{inf}. It has the form
 $$\,
  \rho\,(\rho-1)\,(\rho-2) + 2 \nu\,\rho\,(\rho-1) + \nu\,(\nu-1)\,\rho=0
\,$$ 
and $\rho^{\infty}_1=0,\,\rho^{\infty}_2=1-\nu,\,\rho^{\infty}_3=2-\nu$ are its roots.
	As we mentioned in the introduction
	the characteristic exponents $\rho^{\infty}_i\,,\,i=1, 2, 3$ at $x=\infty$ coincide with the characteristic exponents
	 at the same  point $x=\infty$ of the initial	equation.

	With respect to  the above characteristic exponents $\rho^j_i$ we define
		the following exponent differences:
	  \be\label{ex}
		    	\Delta^L_{32}=\rho^L_3-\rho^L_2= \frac{1}{\sqrt{\varepsilon}}\,,
								  & &
				\Delta^R_{32}=\rho^R_3-\rho^R_2= - \Delta^L_{32}\,,\nonumber\\[0.2ex]
		  \Delta^R_{21}=\rho^R_2-\rho^R_1=\frac{\nu}{2} + \frac{1}{2 \sqrt{\varepsilon}}\,,  
			   & &
			\Delta^L_{21}=\rho^L_2-\rho^L_1=\Delta^R_{21} + \Delta^R_{32}\,,\\[0.2ex]
			\Delta^R_{31}=\rho^R_3-\rho^R_1=\Delta^L_{21}\,,
			  & &
				\Delta^L_{31}=\rho^L_3-\rho^L_1=\Delta^R_{21}\,.\nonumber
		\ee
	
	 Classically, the Fuchsian  singular point $x_j,\,j=R, L$ is called a resonant Fuchsian
		singularity, if there is a $\Delta^j_{kp},\,k \neq p,\,k=2, 3,\,p=1, 2$, which is an integer
		\cite{IKSY}.
		Otherwise, it is called a non-resonant Fuchsian singularity \cite{IKSY}.
		The local theory of the Fuchsian singularities says that the presence of
		a resonant regular singular point $x_j, j=R, L$ is a necessary but not a sufficient condition
		the fundamental matrix $\Phi(x, \varepsilon)$ to contain logarithmic terms \cite{IKSY}.
		On the other hand we have already observed that only the elements $\Phi_{12}(x, \varepsilon)$
		and $\Phi_{13}(x, \varepsilon)$ of the  matrix $\Phi(x, \varepsilon)$ can
		contain logarithmic terms. 
	 So, from here on
		we focus mainly on the computation  of the monodromy matrices $M_j(\varepsilon),\,j=R, L$
		relevant to the so called resonant logarithmic cases (see below).

 For simplicity, through the present and the next section, we
				call these values of the parameters $\nu$ and $\varepsilon$, for which there
	is a chance  of the appearance of logarithmic terms near a resonant Fuchsian singularity $x_j, j=R, L$
	just {\bf the  resonant logarithmic  cases} or
	{\bf a logarithmic resonance}. 
	We distinguish  three  different types of resonant logarithmic cases.
	 They are:
	  \begin{itemize}
	  \item\,
		the resonant logarithmic cases of type ${\bf (B)}$ if
		  $$\,{\bf (B)}\qquad
			 \Delta^R_{21}=\Delta^L_{31}\in\ZZ\,,\qquad \Delta^L_{21}=\Delta^R_{31}\in\ZZ\qquad
			\textrm{simultaneously}\,;
			\,$$
			\item\,
			the  resonant logarithmic cases of type ${\bf (C)}$ if
			$$\,{\bf (C)}\qquad
			 \Delta^L_{21}=\Delta^R_{31}\in\ZZ\,,\qquad \Delta^R_{21}=\Delta^L_{31}\notin\ZZ\qquad
			\textrm{simultaneously}\,;
			\,$$
			\item\,
			the resonant logarithmic cases of type ${\bf (D)}$ if
			$$\,{\bf (D)}\qquad
			 \Delta^R_{21}=\Delta^L_{31}\in\ZZ\,,\qquad \Delta^L_{21}=\Delta^R_{31}\notin\ZZ\qquad
			\textrm{simultaneously}\,.
			\,$$
	\end{itemize}
 It turns out that we can already here reduce the number of the types of the resonant logarithmic cases. 
 Indeed, the simultaneous conditions, which define the resonant logarithmic cases of type
 ${\bf (D)}$, imply that
   $$\,
		  \frac{\nu}{2}+\frac{1}{2 \sqrt{\varepsilon}}\in\ZZ\quad \textrm{but}\quad
			\frac{\nu}{2}-\frac{1}{2 \sqrt{\varepsilon}} \notin \ZZ\,.
		\,$$
  Compare this conditions with the form of the elements $\Phi_{12}(x, \varepsilon)$ and
	$\Phi_{13}(x, \varepsilon)$. We  see that the element $\Phi_{12}(x, \varepsilon)$
	(resp. $\Phi_{13}(x, \varepsilon)$) can contain logarithmic term near $x_R$ (resp. $x_L$)
	if and only if $\nu/2+1/2 \sqrt{\varepsilon}-1 < 0$. On the other hand the limit
	$\sqrt{\varepsilon} \rightarrow 0\in\RR_+$ is equivalent to the limit
	 $\nu/2+1/2 \sqrt{\varepsilon} \rightarrow +\infty$ for a fixed $\nu$. But this implies
	that $\nu/2+1/2 \sqrt{\varepsilon}\in\NN$. Therefore during a logarithmic resonance
	of type ${\bf (D)}$ the fundamental matrix $\Phi(x, \varepsilon)$ does not contain logarithmic terms.
	So, in this article we consider only the rest
		resonant logarithmic cases.
	
		Now, we are in a position to describe the behavior of the
		fundamental matrix $\Phi(x, \varepsilon)$  at the finite singularities
		$x_j,\,j=R, L$ during a logarithmic resonance of type $B$ and $C$.
		We have the following result.
		\bpr{lb}
      During a logarithmic resonance of type ${\bf (B)}$ and ${\bf (C)}$ the fundamental matrix
			$\Phi(x, \varepsilon)$ of the perturbed equation  is represented
			near the singular points $x_j,\,j=R, L$, as
			  \ben
				 \Phi(x, \varepsilon)=
					 (I_L(\varepsilon) + \mathcal{O}(x-x_L))\,
					(x-x_L)^{\frac{1}{2} \Lambda + \frac{1}{2 x_L} Q}\,(x-x_L)^{T_L}
				\een
				in a neighborhood of $x_L$, which does not contain the point $x_R$, and
					\ben
					 \Phi(x, \varepsilon)=
					 (I_R(\varepsilon) + \mathcal{O}(x-x_R))\,
					(x-x_R)^{\frac{1}{2} \Lambda + \frac{1}{2 x_R} Q}\,(x-x_R)^{T_R}
					\een
					in a neighborhood of $x_R$, which does not contain the point $x_L$.
			The matrices $I_j(\varepsilon) + \mathcal{O}(x-x_j)$ are analytic matrix-functions
			near the point $x_j,\,j=R, L$, respectively.
			The matrices $\Lambda$ and  $Q$  are the matrices, associated with the initial equation
			and defined by \eqref{l},
			  \be\label{J}
				T_j=\left(\begin{array}{ccc}
										0     &d^j_2      &d^j_3\\
										0     &0          &0\\
										0     &0          &0
									\end{array}
						\right)\,.				
				\ee
			The elements $d^j_i,\,i=2, 3,\,j=R, L$ of the matrices $T_j,\,j=R, L$
			are defined as follows,
			  					 \be\label{d}
					  d^R_2=0\,,& &  \quad 
						d^L_2=\res \left( 
						 \frac{(x-\sqrt{\varepsilon})^{\frac{1}{2 \sqrt{\varepsilon}}+\frac{\nu-2}{2}}}
				     {(x+\sqrt{\varepsilon})^{\frac{1}{2 \sqrt{\varepsilon}}-\frac{\nu-2}{2}}}\,,\,x=x_L\right)\,,\\[0.5ex]
						   d^L_3=0\,,& & \quad	  
						d_3^R=-\frac{1}{2}\res\left(
			 \frac{(x+\sqrt{\varepsilon})^{\frac{1}{2 \sqrt{\varepsilon}}+\frac{\nu-2}{2}}}
				     {(x-\sqrt{\varepsilon})^{\frac{1}{2 \sqrt{\varepsilon}}-\frac{\nu-2}{2}}}\,,\,
						x=x_R\right)\,.\nonumber
						\ee
				\epr
	
    \proof
		
		  To study the behavior of the fundamental matrix $\Phi(x, \varepsilon)$
     we will use its explicit form given by \thref{pe-fms}.
		 Let us firstly consider the elements $\Phi_s(x, \varepsilon),\,
		s=2, 3, 23$. We have 
		  \ben
			  \Phi_s(x, \varepsilon)=\left\{\begin{array}{cl}
				  (x-x_L)^{m_s^L}\,g_{s, L}(x)   &\textrm{near the point}\,\,\,x_L\,,\\[0.15ex]
					(x-x_R)^{m_s^R}\,g_{s, R}
					&\textrm{near the point}\,\,\, x_R\,,
					                            \end{array}
																\right.			
			\een
		 where $s=2, 3, 23$, the exponents $m_2^j=\rho_2^j-1,\, m_3^j=m_{23}^j=\rho_3^j-2,\,j=L, R$.
		The functions $g_{s, j}(x)$ are analytic functions in $x$ in a neighborhood of the
		singular point $x_j$, which does not contain the other singular point.
		In particular $g_{s, j}(x_j) \neq 0$ for $s=2, 3$, and $g_{23, j}(x_j)=0$.
			
			Now we consider the elements $\Phi_s(x, \varepsilon)$ for $s=1, 12, 13$.
			We have
			\ben
				 \Phi_1(x, \varepsilon)=\left\{\begin{array}{cl}
				     (x-x_L)^{\rho_1^L}\,g_{1, L}(x)   &\textrm{near the point}\,\,\,x_L\,,\\[0.15ex]
						(x-x_R)^{\rho^R_1}\,g_{1, R}(x)    &\textrm{near the point}\,\,\,x_R\,.
						                          \end{array}
																\right.			
				\een
				For the element $\Phi_{12}(x, \epsilon)$ we have
			   \ben
				   \Phi_{12}(x, \varepsilon)=
	         d^L_2\,(x-x_L)^{\rho^L_1}\,\log(x-x_L)\,g_{1, L}(x) + (x-x_L)^{\rho^L_2-1}\,h_{2, L}(x)
					\een
					in a neighborhood of $x_L$, which does not contain $x_R$, and
					\ben
					\Phi_{12}(x, \varepsilon)=			
										(x-x_R)^{\rho^R_2-1}\,h_{2, R}(x)
									\een
						in a neighborhood of $x_R$, which does not contain the point $x_L$.			
		  For the element $\Phi_{13}(x, \varepsilon)$ we have
			 \ben
			  \Phi_{13}(x, \varepsilon) =
				(x-x_L)^{\rho^L_3-2}\,h_{3, L}(x)
				\een
				near the point $x_L$, and
				\ben
				 \Phi_{13}(x, \varepsilon)=
					  d^R_3\,(x-x_R)^{\rho^R_1}\,\log(x-x_R)\,g_{1, R}(x) + 
					(x-x_R)^{\rho^R_3-2}\,h_{3, R}(x)
				\een
				in a neighborhood of $x_R$, which does not contain the point $x_L$. 
			The functions $g_{1, j}(x)$ and $h_{s, j},\,j=L, R,\,s=2, 3$ are analytic functions in $x$ 
			in a neighborhood of the singular point $x_j$, which does not contain the other singular point.
			In particular, $g_{1, j}(x_j) \neq 0$, and $h_{s, j}(x_j)=0$.
			
		Note that the absence of logarithmic terms during a resonance implies nullity
		of the numbers  $d^j_i,\,1=2, 3$.\qed

		 \bre{solution1}
		From \prref{lb} it follows that, the norm of the $j$-column of the
		fundamental matrix $\Phi(x, \varepsilon)$ has the same growth rate near the singular points
		as the function $\Phi_j(x, \varepsilon),\,j=1, 2, 3$.
	  In particular, near $x_R$ the function $\Phi_2(x, \varepsilon)$ 
		has the smallest growth, and the function $\Phi_3(x, \varepsilon)$ has the largest one.
		Near $x_L$ the function $\Phi_3(x, \varepsilon)$ has the smallest growth, and the
		function $\Phi_2(x, \varepsilon)$ has the largest one. The first column is determined
		as the unique solution (up to a constant factor) that has a mid-growth at both $x_R$
		and $x_L$.
				\ere

   \bth{monodromy}
	   During a logarithmic resonance  of type ${\bf (B)}$ and ${\bf (C)}$
		the monodromy matrices $M_j(\varepsilon),\,j=R, L$ of the perturbed equation
		with respect to the fundamental matrices $\Phi(x, \varepsilon)$,
		defined by \thref{pe-fms}, are given by
				\be\label{B-C-D}
					  M_j(\varepsilon)=e^{\pi\,i\,(\Lambda+\frac{1}{x_j} Q)}\,e^{2 \pi\,i\,T_j}=
						e^{2 \pi\,i\,T_j}\,e^{\pi\,i\,(\Lambda+\frac{1}{x_j} Q)}\,. 
					\ee	
	\ethe

     \proof
		
		  The proof follows from \prref{lb} and \deref{mr}, \deref{monodromy} together
			with the observation that during a resonance the matrices 
			$e^{\pi\,i\,(\Lambda+\frac{1}{x_j} Q)}$ and $(x-x_j)^{T_j}$ commute. 
			Note that during the logarithmic resonance of type ${\bf (B)}$ and ${\bf (C)}$
			the matrices $e^{2 \pi\,i\,T_j}$ and $e^{\pi\,i\,(\Lambda+\frac{1}{x_j} Q)}$ also
			commute. \qed

		\bre{monodromy}
		   Combining \prref{lb} and \thref{monodromy} we observe 
			that the first two columns of the matrix $\Phi(x, \varepsilon)$
			are eigenvectors of the monodromy operators $M_R(\varepsilon)$ with eigenvalues 
			$e^{2 \pi\,i\,\rho^R_1}=e^{\pi\,i/\sqrt{\varepsilon}}$ and
			$e^{2 \pi\,i\,(\rho^R_2-1)}=e^{\pi\,i(\nu+1/\sqrt{\varepsilon})}$, respectively.
			In the same manner the first column and the third column of the matrix $\Phi(x, \varepsilon)$  
			are eigenvectos of the monodromy operator $M_L(\varepsilon)$ with eigenvalues
			$e^{2 \pi\,i\,\rho^L_1}=e^{-\pi\,i/\sqrt{\varepsilon}}$ and
			$e^{2 \pi\,i\,(\rho^L_2-1)}=e^{\pi\,i(\nu-1/\sqrt{\varepsilon})}$, respectively.
			The numbers $d^R_3$ and $d^L_2$ (when they are different from zero) block
			 the third  and the second column of $\Phi(x, \varepsilon)$
			 of being eigenvectors of the monodromy operators
			$M_R(\varepsilon)$ and $M_L(\varepsilon)$ respectively (see also \cite{CL-CR2, Z} about the confluence
			of the hypergeometric equation).
			
				Combining the observations of \reref{solution1} and \reref{monodromy}, we see that
				the columns of the fundamental matrix  of the perturbed equation leads to the so called
				``mixed basis'' of solutions (see \cite{HLR}, Theorem 5.4 and 5.5 for more details).
		\ere
In what follows we will calculate explicitly the numbers
$d^j_i,\,i=2, 3,\,j=R, L$ \eqref{d}	 and write down the corresponding
monodromy matrices \eqref{B-C-D}.


		\subsection{The resonant logarithmic cases of type ${\bf (B)}$.}

   We begin our calculations with  the resonant logarithmic cases of type  ${\bf (B)}$.

 The simultaneous conditions 
	$\Delta^R_{21}=\Delta^L_{31}\in\ZZ,\,\Delta^L_{21}=\Delta^R_{31}\in\ZZ$ imply that
	 \ben
	   \nu\in\ZZ,\quad \frac{1}{\sqrt{\varepsilon}}\in\NN\quad
		 \textrm{such that}\quad
		 \frac{\nu}{2}+\frac{1}{2 \sqrt{\varepsilon}}\in\ZZ\,.
	 \een
	 In the last section we are going to study the behavior of
	the monodromy matrices, obtained in the present section, when 
	$\sqrt{\varepsilon} \rightarrow 0\in\RR_+$ and $1/\sqrt{\varepsilon}\in\NN$ for fixed $\nu$. But the limit 
	when $\sqrt{\varepsilon} \rightarrow 0\in\RR_+$
	(for fixed $\nu$) is equivalent to the limit $\nu/2+1/2 \sqrt{\varepsilon} \rightarrow +\infty$. So,
	we limit calculation to the case when $\nu/2+1/2 \sqrt{\varepsilon}\in\NN$. 
	We also note that when
	$\nu/2 +1/2 \sqrt{\varepsilon}\in\NN$ the integral 
	$\Phi_{12}(x, \varepsilon)/\Phi_1(x, \varepsilon)$ in \eqref{lsol1}
	is a convergent one.
	
\bth{B}
  Assume that  $\nu\in\ZZ$. Assume also that  
	$1/\sqrt{\varepsilon}\in\NN$ 
	such that $\nu/2+1/2 \sqrt{\varepsilon}\in\NN$. 
	Then  for $d^j_i,\,i=2, 3,\,j=R, L$ \eqref{d} and the corresponding monodromy matrices $M_j$,  $j=R, L$ 
	\eqref{B-C-D} we  have 
	   \ben
	   d^R_2=d^L_3=0\,,& & \quad
		 d^L_2=\left(-\frac{1}{2 \sqrt{\varepsilon}}\right)^{1-\nu}
		\frac{\Gamma(\frac{1}{2 \sqrt{\varepsilon}}+\frac{\nu}{2})}
		     {\Gamma(\nu)\,\Gamma(\frac{1}{2 \sqrt{\varepsilon}} - \frac{\nu}{2}+1)}\,,\\[0.4ex]
				          & &
				d^R_3=-\frac{1}{2}\left(\frac{1}{2 \sqrt{\varepsilon}}\right)^{1-\nu}
		\frac{\Gamma(\frac{1}{2 \sqrt{\varepsilon}}+\frac{\nu}{2})}
		     {\Gamma(\nu)\,\Gamma(\frac{1}{2 \sqrt{\varepsilon}} - \frac{\nu}{2}+1)}\,. 
		 \een
		The corresponding monodromy matrices  are given by
		 $$\,
			   M_j(\varepsilon)=\left(\begin{array}{ccc}
				  (-1)^{\nu}     &2 \pi i (-1)^{\nu}\,d^j_2             &2 \pi i (-1)^{\nu}\,d^j_3\\
					0              &(-1)^{\nu}     &0\\
					0              &0              &(-1)^{\nu}
					        \end{array}
						\right)\,,\quad j=R, L\,.
		\,$$				
				\ethe

\proof

According to \eqref{d} the numbers  $d^j_2,\, j=R, L$ are defined by 
  $$\,
	  d^j_2=\res \left(\frac{(x-\sqrt{\varepsilon})^{\frac{1}{2 \sqrt{\varepsilon}}+\frac{\nu}{2}-1}}
		{(x+\sqrt{\varepsilon})^{\frac{1}{2 \sqrt{\varepsilon}}-\frac{\nu}{2}+1}}\,,\,x=x_j\right)\,,\qquad
		j=R, L\,.
	\,$$
	Since $1/2 \sqrt{\varepsilon}+\nu/2 \in\NN$ we find that $d^R_2=0$.
	When $\nu\in\ZZ$ and $1/2 \sqrt{\varepsilon}+\nu/2\in\NN$ the exponent
	$1/2\sqrt{\varepsilon}-\nu/2+1\in\ZZ$. Then for the number $d^L_2$ we ontain 
	consecutively
	 \ben
	   d^L_2   &=&
		\frac{1}{\left(\frac{1}{2 \sqrt{\varepsilon}}-\frac{\nu}{2}\right)!}
		\left(
		   \frac{d^{\frac{1}{2 \sqrt{\varepsilon}}-\frac{\nu}{2}}}
			{d x^{\frac{1}{2 \sqrt{\varepsilon}}-\frac{\nu}{2}}} 
			(x-\sqrt{\varepsilon})^{\frac{1}{2 \sqrt{\varepsilon}}+\frac{\nu}{2}-1}
			                     \right)_{x=x_L}=\\[0.4ex]
						&=&							
				\left(-\frac{1}{2 \sqrt{\varepsilon}}\right)^{1-\nu}
		\frac{\Gamma(\frac{1}{2 \sqrt{\varepsilon}}+\frac{\nu}{2})}
		     {\Gamma(\nu)\,\Gamma(\frac{1}{2 \sqrt{\varepsilon}} - \frac{\nu}{2}+1)}\,,
	\een
	 where $\Gamma(z)$ is the classical Gamma function.
	
	Next, we compute the numbers $d^j_3$. Applying the formula \eqref{d}, we see that
		$d^L_3=0$ when $1/2 \sqrt{\varepsilon}+\nu/2\in\NN$. For the number $d^R_3$
		we have
  \ben
	  d^R_3    &=&											
					-\frac{1}{2\left(\frac{1}{2 \sqrt{\varepsilon}}-\frac{\nu}{2}\right)!}
					   \left(\frac{d^{\frac{1}{2 \sqrt{\varepsilon}}-\frac{\nu}{2}}}
			{d x^{\frac{1}{2 \sqrt{\varepsilon}}-\frac{\nu}{2}}} 
			(x-\sqrt{\varepsilon})^{\frac{1}{2 \sqrt{\varepsilon}}+\frac{\nu}{2}-1}
			                     \right)_{x=x_L}=\\[0.4ex]
						&=&							
				-\frac{1}{2}\left(\frac{1}{2 \sqrt{\varepsilon}}\right)^{1-\nu}
		\frac{\Gamma(\frac{1}{2 \sqrt{\varepsilon}}+\frac{\nu}{2})}
		     {\Gamma(\nu)\,\Gamma(\frac{1}{2 \sqrt{\varepsilon}} - \frac{\nu}{2}+1)}\,.
									\een
	We also note that in this case  the integral, including in the definition of $\Phi_{13}(x, \varepsilon)$
	 $$\,
	  \int_{-\sqrt{\varepsilon}}^x
		     \frac{(t+\sqrt{\varepsilon})^{\frac{1}{2 \sqrt{\varepsilon}}+\frac{\nu}{2} -1}}
					    {(t-\sqrt{\varepsilon})^{\frac{1}{2 \sqrt{\varepsilon}} -\frac{\nu}{2}+1}}\, d t
	\,$$
	is a convergent one.

 Applying \thref{monodromy}, we see that the corresponding monodromy
are given by
	$$\,
	   M_j(\varepsilon)=\left(\begin{array}{ccc}
				  (-1)^{\nu}     &2 \pi i (-1)^{\nu}\,d^j_2             &2 \pi i (-1)^{\nu}\,d^j_3\\
					0              &(-1)^{\nu}     &0\\
					0              &0              &(-1)^{\nu}
					        \end{array}
						\right)\,,\quad j=R, L\,,
		\,$$				
since $\nu$ is an integer.

This completes the proof of \thref{B}.\qed

\vspace{1ex}


\subsection{The resonant logarithmic cases of type ${\bf (C)}$ .}

The conditions $\Delta^L_{21}=\Delta^R_{31}\in\ZZ$ but $\Delta^R_{21}=\Delta^L_{31}\notin\ZZ$
which define the resonant logarithmic case ${\bf (C)}$ imply that
 $$\,
  \frac{\nu}{2} - \frac{1}{2 \sqrt{\varepsilon}}\in\ZZ\quad
	\textrm{but} \quad \nu\notin\ZZ\,.
\,$$
Similar to the resonant logarithmic cases of type ${\bf B}$, we can restrict the
calculations to the non-positive values of $\nu/2-1/2 \sqrt{\varepsilon}$.

Then we have

\bth{C}
  Assume that $1/2 \sqrt{\varepsilon}-\nu/2\in\NN$  but $\nu\notin \ZZ$.
	Then for the numbers $d^j_i$ and the corresponding  monodromy matrices $M_R$ and $M_L$ 
of the perturbed equation  we have
		 \ben
	   d^R_2=d^L_3=0\,,& & \quad
		 d^L_2=\left(-\frac{1}{2 \sqrt{\varepsilon}}\right)^{1-\nu}
		\frac{\Gamma(\frac{1}{2 \sqrt{\varepsilon}}+\frac{\nu}{2})}
		     {\Gamma(\nu)\,\Gamma(\frac{1}{2 \sqrt{\varepsilon}} - \frac{\nu}{2}+1)}\,,\\[0.4ex]
				          & &
				d^R_3=-\frac{1}{2}\left(\frac{1}{2 \sqrt{\varepsilon}}\right)^{1-\nu}
		\frac{\Gamma(\frac{1}{2 \sqrt{\varepsilon}}+\frac{\nu}{2})}
		     {\Gamma(\nu)\,\Gamma(\frac{1}{2 \sqrt{\varepsilon}} - \frac{\nu}{2}+1)}\,. 
		 \een
		  The monodromy matrices are given by
		  \ben
			  M_R(\varepsilon)=\left(\begin{array}{ccc}
				  e^{\pi\,i \nu}    &0                   
					&2 \pi\,i e^{ \pi\,i \nu} d^R_3\\
					 0                 &e^{3 \pi\,i \nu}   &0\\
					 0                 &0                   &e^{\pi\,i \nu}
					        \end{array}
						\right)\,,
			M_L(\varepsilon)=\left(\begin{array}{ccc}
				  e^{-\pi\,i \nu}    &2 \pi\,i e^{-\pi\,i \nu} d^L_2                   &0\\
					 0                 &e^{-\pi\,i \nu}   &0\\
					 0                 &0                   &e^{\pi\,i \nu}
					        \end{array}
						\right)\,.						
			\een
\ethe

\proof

As in the previous theorems for the numbers $d^j_3$ we have that
 \ben
  d^j_3 =-\frac{1}{2}
	\res\left( 
	\frac{(x+\sqrt{\varepsilon})^{\frac{1}{2 \sqrt{\varepsilon}}+\frac{\nu}{2} -1}}
	{(x-\sqrt{\varepsilon})^{\frac{1}{2\sqrt{\varepsilon}} -\frac{\nu}{2}+1}}\,,\,
	x=x_j\right)\,.
	\een
Since the exponent $1/2 \sqrt{\varepsilon}+\nu/2-1 \notin \ZZ$ when $\nu\notin\ZZ$ then $d^L_3=0$. 
We also note that for sufficiently big $1/\sqrt{\varepsilon}$ 
 the integral, including in the definition of $\Phi_{13}(x, \varepsilon)$
$$\,
  \int_{-\sqrt{\varepsilon}}^x 
	\frac{(t+\sqrt{\varepsilon})^{\frac{1}{2 \sqrt{\varepsilon}}+\frac{\nu}{2} -1}}
	{(t-\sqrt{\varepsilon})^{\frac{1}{2\sqrt{\varepsilon}} -\frac{\nu}{2}+1}}\,d t
	\,$$
is a convergent one. Next, for the number $d^R_3$ we find
  $$\,
				d^R_3=-\frac{1}{2}\left(\frac{1}{2 \sqrt{\varepsilon}}\right)^{1-\nu}
		\frac{\Gamma(\frac{1}{2 \sqrt{\varepsilon}}+\frac{\nu}{2})}
		     {\Gamma(\nu)\,\Gamma(\frac{1}{2 \sqrt{\varepsilon}} - \frac{\nu}{2}+1)}\,. 
	\,$$

	For the number $d^L_2$ we obtain
 $$\,
 d^L_2=\left(-\frac{1}{2 \sqrt{\varepsilon}}\right)^{1-\nu}
		\frac{\Gamma(\frac{1}{2 \sqrt{\varepsilon}}+\frac{\nu}{2})}
		     {\Gamma(\nu)\,\Gamma(\frac{1}{2 \sqrt{\varepsilon}} - \frac{\nu}{2}+1)}\,.
 \,$$
Again since the exponent $1/2 \sqrt{\varepsilon}+\nu/2-1 \notin \ZZ$, then $d^R_2=0$.
This ends the calculations of the numbers $d^j_i$.

Finally, applying \eqref{B-C-D} and using the connection
$1/2 \sqrt{\varepsilon}-\nu/2\in\NN$, we write down the
monodromy matrices $M_R(\varepsilon)$ and $M_L(\varepsilon)$.
\qed

\bre{E2}
  Similar to the 	non-perturbed equation when $\nu=0$ the point $x=\infty$ is 
	no longer a singular point for the 	perturbed equation. Indeed, when $\nu=0$ we
	set $x=1/t$. Then the 	perturbed equation becomes
	 \ben
	    & &
	  \dddot{y}(t) +\left[
		-3 \sqrt{\varepsilon}\left(1 - \frac{1}{2 \sqrt{\varepsilon}}\right)\,\frac{1}{1-\sqrt{\varepsilon} t} +
		3 \sqrt{\varepsilon}\left(1+\frac{1}{2 \sqrt{\varepsilon}}\right)\,\frac{1}{1+\sqrt{\varepsilon} t}
		              \right]\,\ddot{y}(t) +\\[0.4ex]
			&+&
	\left[\frac{-2 \varepsilon + \frac{1}{2}}{1-\sqrt{\varepsilon} t} +
				\frac{-2 \varepsilon + \frac{1}{2}}{1+\sqrt{\varepsilon} t} +
				\frac{\varepsilon (1-\frac{3}{2 \sqrt{\varepsilon}} + \frac{1}{2 \varepsilon})}{(1-\sqrt{\varepsilon})^2}+
				\frac{\varepsilon (1+\frac{3}{2 \sqrt{\varepsilon}} + \frac{1}{2 \varepsilon})}{(1+\sqrt{\varepsilon})^2}
	\right]\,\dot{y}(t)=0								
	\een
	for which the point $t=0$ (resp. the point $x=\infty$ for the original equation) is an 
	ordinary point.
\ere

  
	\section{ Main results }
	
	In this section we will connect, by a limit $\sqrt{\varepsilon} \rightarrow 0$,
   analytic invariants of the initial and the perturbed equations,
	computed in the previous two sections. 
	In order to connect by a  limit $\sqrt{\varepsilon} \rightarrow 0$
		the solution of the perturbed equation with the solution of the initial equation,
		as well as their invariants, we consider the perturbed equation on the sectorial domains
		$\Omega_1(\varepsilon)$ and $\Omega_2(\varepsilon)$. These domains are obtained from
		$\Omega_1$ and $\Omega_2$ (relevant to the initial equation) by making a cut
		between the singular points $x_L$ and $x_R$ through the real axis, The point $x_0=0$
		belongs to this cut. When $\varepsilon \rightarrow 0$ then $\Omega_j(\varepsilon), j=1, 2$
		tend to $\Omega_j, j=1, 2$.
		The domains $\Omega_1(\varepsilon)$ and $\Omega_2(\varepsilon)$ intersect in the left $\Omega_L(\varepsilon)$, 
		right $\Omega_R(\varepsilon)$ sectors and along the cut. The points $x_j, j=L, R$ do
		not belong to $\Omega_j(\varepsilon)$, but belong to their closure, respectively. Both poits
		belong to the cut (see \cite{CL-CR}).

	  In the keeping with the initial equation, we rewrite the fundamental matrix
	$\Phi(x, \varepsilon)$ of the perturbed equation in the form 
	 $$\,
	   \Phi(x, \varepsilon)=H(x, \varepsilon)\,F(x, \varepsilon)\,.
	\,$$
	 Here $F(x, \varepsilon)=(x-x_L)^{\Lambda/2 + Q/2 x_L}\,(x-x_R)^{\Lambda/2+Q/2 x_R}$.
	 The matrix $H(x, \varepsilon)$ is defined as
	  \ben
		  H(x, \varepsilon)=\left(\begin{array}{ccc}
			1   &\frac{\Phi_1(x, \varepsilon)}{\Phi_2(x, \varepsilon)}
			  \int_{\Gamma_1(x, \varepsilon)}\frac{\Phi_2(t, \varepsilon)}{\Phi_1(t, \varepsilon)}\,d t
				 &\frac{\Phi_1(x, \varepsilon)}{\Phi_3(x, \varepsilon)}
			  \int_{\Gamma_2(x, \varepsilon)}\frac{\Phi_{23}(t, \varepsilon)}{\Phi_1(t, \varepsilon)}\,d t\\[0.5ex]
		0   &1  &-\frac{x^2-\varepsilon}{2}\\[0.3ex]
		0   &0  &1
				                 \end{array}
												\right)\,.
			\een												
		When we move the path $\Gamma_1(x, \varepsilon)$ analytically around the origin we obtain near
		the negative real axis $\RR_-$ two branches of the element $h_{12}(x, \varepsilon)$ 
		of the matrix $H(x, \varepsilon)$. We denote them by		
		$h_{12}^+(x, \varepsilon)$ and $h_{12}^-(x, \varepsilon)$. The function $h_{12}^-(x, \varepsilon)$
		is defined on the path in the direction $\pi-\epsilon$, and the function $h_{12}^+(x, \varepsilon)$
		is defined in the direction $\pi+\epsilon$. 
		Similarly, near the positive real axis $\RR_+$ we have two branches, 
		$h_{13}^+(x, \varepsilon)$ and $h_{13}^-(x, \varepsilon)$, of the 
		function $h_{13}(x, \varepsilon)$. When $\Gamma_1(x, \varepsilon)$
		(resp. $\Gamma_2(x, \varepsilon)$) crosses $\RR_-$ (resp. $\RR_+$) we rather observe 
		a Stokes phenomenon, that a linear monodromy. This phenomenon is described by the so called
		unfolded Stokes matrices. In accordance with the initial equation, we
		determine on the sector $\Omega_1(x, \varepsilon)$ the fundamental matrix of the perturbed equation as
		$$\,
		  \Phi_1(x, \varepsilon)=H_1(x, \varepsilon)\,F_1(x, \varepsilon)\,,
		\,$$
		 where
		 $$\,
		    H_1(x, \varepsilon)=\left(\begin{array}{ccc}
				  1   &h^-_{12}(x, \varepsilon)   &h^+_{13}(x, \varepsilon)\\
					0   &1                          &h_{23}(x, \varepsilon)\\
					0   &0                          &1
					                         \end{array}
														\right)\,,			
		\,$$
		and $F_1(x, \varepsilon)$ is the branch of $F(x, \varepsilon)$ on $\Omega_1(\varepsilon)$.
		On the sector $\Omega_2(\varepsilon)$ we define the fundamental matrix as
		 $$\,
		  \Phi_2(x, \varepsilon)=H_2(x, \varepsilon)\,F_2(x, \varepsilon)\,,
		\,$$
		where
		  $$\,
		    H_2(x, \varepsilon)=\left(\begin{array}{ccc}
				  1   &h^+_{12}(x, \varepsilon)   &h^-_{13}(x, \varepsilon)\\
					0   &1                          &h_{23}(x, \varepsilon)\\
					0   &0                          &1
					                         \end{array}
														\right)\,,
		\,$$
	as $F_2(x, \varepsilon)=F_1(x, \varepsilon)$ on $\Omega_L(\varepsilon)$ 
	and $F_2(x, \varepsilon)=F_1(x, \varepsilon)\,\hat{M}$ on $\Omega_R(\varepsilon)$.
	  Denote by $St_L(\varepsilon)$ and $St_R(\varepsilon)$ the unfolded Stokes matrices with respect
		to the fundamental matrix $\Phi_1(x, \varepsilon)$ on the upper sector $\Omega_1(\varepsilon)$.
		Now we will describe the change of the fundamental matrix when we turn around the origin
		analytically in the positive sense. We start from the sector $\Omega_1(\varepsilon)$
			and the solution $\Phi_1(x, \varepsilon)$ on it. When $\Gamma_1(x, \varepsilon)$ crosses
		the negative real axis, we observe a Stokes phenomenon on $\Omega_L(\varepsilon)$. In particular, 	
		the unfolded Stokes matrix $St_L(\varepsilon)$ is defined by 
		 \ben
		  St_L(\varepsilon)=(\Phi_2(x, \varepsilon))^{-1}\,\Phi_1(x, \varepsilon) \quad
				\textrm{on}\quad \Omega_L(\varepsilon)
		\een
		If we continue circling round the origin then when $\Gamma_2(x, \varepsilon)$ crosses
		the positive real axis, we  observe  a Stokes phenomenon on $\Omega_R(\varepsilon)$.
		 The jump of the solution
		$\Phi_2(x, \varepsilon)$  to the solution $\Phi_1(x, \varepsilon)$ is defined by
		 \ben
		 (\Phi_1(x, \varepsilon))^{-1}\,\Phi_2(x, \varepsilon)=
			   St_R(\varepsilon)\,\hat{M} \quad
					\textrm{on}\quad \Omega_R(\varepsilon)\,,
		\een
		since on $\Omega_2(\varepsilon)$ we have $F_2(x, \varepsilon)=F_1(x, \varepsilon)\,\hat{M}$.

	From Theorem 4.25 in \cite{CL-CR} it follows that the unfolded Stokes matrices
	$St_L(\varepsilon)$ and $St_R(\varepsilon)$ depend analytically on the parameter of
		perturbation $\varepsilon$ and they converge when $\varepsilon \rightarrow 0$ to
		the Stokes matrices $St_L=St_{\pi}$ and $St_R=St_0$ of the initial equation.

		In the previous section we have computed the monodromy matrix $M_R(\varepsilon)$ of
		the perturbed equation with respect to the fundamental solution, defined on the upper
		sector $\Omega_1(\varepsilon)$. 
		With respect to the fundamental solution on the lower sector $\Omega_2(\varepsilon)$
		the  monodromy matrix $\tilde{M}_R(\varepsilon)$ is given by 
	$$\,
		 \tilde{M}_R(\varepsilon)=\hat{M}^{-1}\,M_R(\varepsilon)\,\hat{M}\,,
		\,$$
		where $M_R(\varepsilon)$ is the monodromy matrix, defined by \eqref{B-C-D}.  
		
		Now, we can give the connection between the monodromy matrices 
		and the unfolded Stokes matrices. Proposition 4.31 in \cite{CL-CR} states that the monodromy
		operator acting on the solution $\Phi_j(x, \varepsilon)$ decomposes into the Stokes
		operator  multiplied, from the right, by the
		classical monodromy operator acting on branch of $F(x, \varepsilon)$. 
    In \cite{MK3} Theorem 32, Klime\v s expresses in a remarkable way the acting of the monodromy
		operators on analytic extension of the solutions of the perturbed equation to the whole
		$\Omega_1(\varepsilon) \cup \Omega_2(\varepsilon)$ by the monodromy matrices $M_j(\varepsilon)$,
		unfolded Stokes matrices $St_j(\varepsilon)$ and the matrices $e^{\pi\,i(\Lambda+Q/x_j)},\,
		j=L, R$. His formulas have been deduced provided that there is an agreement of the matrices
		$F(x)$ and $F(x, \varepsilon)$ on the right intersections $\Omega_R$ and $\Omega_R(\varepsilon)$.
		In the next proposition we  reformulate (without giving a proof) his formulas, provided that
			the above agreement is on the left intersections (see for details and proof \cite{MK3}).

		  \bpr{S-M}
			  Let $M_j(\varepsilon)$ and $St_j(\varepsilon),\, j=R, L$ be the monodromy matrices and 
				 the infolded Stokes matrices of the perturbed equation with respect to the fundamental
				solution on the upper sector $\Omega_1(\varepsilon)$.
				Then on the upper sector $\Omega_1(\varepsilon)$ they satisfy the following relations
				 \ben
				   M_L(\varepsilon)=
					e^{\pi\,i(\Lambda + \frac{1}{x_L} Q)}\,St_L(\varepsilon),\qquad
					 M_R(\varepsilon)=
					St_R(\varepsilon)\,e^{\pi\,i(\Lambda+ \frac{1}{x_R} Q)}\,.
					\een
					 On the lower sector $\Omega_2(\varepsilon)$ they satisfy the following relations
				 \ben
				   M_L(\varepsilon)=
					St_L(\varepsilon)\,e^{\pi\,i(\Lambda + \frac{1}{x_L} Q)},\qquad
					 M_R(\varepsilon)=
					e^{\pi\,i(\Lambda+\frac{1}{x_R} Q)}\,St_R(\varepsilon)\,.
					\een
			\epr
        Note that these relations are in concordance with the definition of the monodromy
			around $x=\infty$ for both equations. Indeed, from \prref{S-M} it follows that
	    on the lower sector $\Omega_2(\varepsilon)$
			 \ben
			   St_L(\varepsilon)\,St_R(\varepsilon)\,\hat{M} &=&
				M_L(\varepsilon)\,\hat{M}^{-1}\,M_R(\varepsilon)\,\hat{M}=M_L(\varepsilon)\,\tilde{M}_R(\varepsilon)=
				M^{-1}_{\infty}(\varepsilon)=\\[0.15ex]
				                &=&
					e^{2 \pi\,i\,T_L}\,e^{2 \pi\,i\,T_R}\,\hat{M}\,.
			\een
    When $\sqrt{\varepsilon} \rightarrow 0$ the monodromy $M^{-1}_{\infty}(\varepsilon)$ around $x=\infty$ of
		the initial equation tends to $St_{\pi}\,St_0\,\hat{M}$. Since the matrices $T_j$ are convergent under this limit,
		then the monodromy matrices $M^{-1}_{\infty}$ of the perturbed equation is well defined. Recall that
		the exponents $\rho^{\infty}_i$ at $x=\infty$ do not changes under the perturbation. So the above
		phenomenon is expected. It is interesting to study if this connection remains valid under perturbation,
		that makes the characteristic exponents $\rho^{\infty}_i$ at $x=\infty$ dependent on the parameter of perturbation.

		\bre{k}
			In fact Theorem 32 in \cite{MK3} states that on the lower sector $\Omega_2(\varepsilon)$
			the monodromy matrix $\tilde{M}_R(\varepsilon)$ is expressed as
				$$\,
					 \tilde{M}_R(\varepsilon)=e^{-\pi\,i(\Lambda+\frac{1}{x_L} Q)}\,St_R(\varepsilon)\,\hat{M}.
				\,$$	
			Using the relation $\tilde{M}_R(\varepsilon)=\hat{M}^{-1}\,M_R(\varepsilon)\,\hat{M}$, we
				rewrite it as $M_R(\varepsilon)=e^{\pi\,i (\Lambda+Q/x_R)}\,St_R(\varepsilon)$. 
		\ere

    In a consequence of \prref{S-M} and \thref{monodromy} we have the following relation.
				\bpr{S-T}
					The unfolded Stokes matrices $St_j(\varepsilon)$ and the matrices $e^{2 \pi\,i\,T_j},\,j=L, R$
					satisfy the following relation
					\ben
					  St_L(\varepsilon)=e^{2 \pi\,i\,T_L}\,,\quad
						St_R(\varepsilon)=e^{2 \pi\,i\,T_R}\,.
					\een
        \epr
				
				\proof
				  From \prref{S-M} we have that
					  $$\,
					 M_R(\varepsilon)=
					St_R(\varepsilon)\,e^{\pi\,i(\Lambda+ \frac{1}{x_R} Q)}\,.
						\,$$
						On other hand the monodromy matrix $M_R(\varepsilon)$, given by \eqref{B-C-D}, is
						\ben
						  M_R(\varepsilon)=
							e^{2 \pi\,i\,T_R}\,e^{\pi\,i(\Lambda+\frac{1}{x_R} Q)}\,.
							\een
						Then comparing the both expressions for $M_R(\varepsilon)$ we obtain the 
						relation between $St_R(\varepsilon)$ and $e^{2 \pi \,i\,T_R}$.
						In the same manner one have the direct relation between $St_L(\varepsilon)$ and
						$e^{2 \pi\,i\,T_L}$.
				\qed

    It turns out that the matrices $T_j, j=L, R$ convergent when $\sqrt{\varepsilon} \rightarrow 0\in\RR_+$.		
   The next  preliminary lemma deal with the limits of the numbers  $d^L_3$ and $d^L_2$, obtained 
	in \thref{B} and \thref{C}.
	
   \ble{first}
			Assume that $\nu\in\RR$ is fixed. 
			Then the numbers $d^L_2$  and $d^R_3$ derived in \thref{B} and \thref{C} satisfy the
			following limits
			  	\ben
		  \lim_{1/\sqrt{\varepsilon}\rightarrow + \infty} d^R_3=
			   -\frac{1}{2 \Gamma(\nu)}\,,\quad				
			\lim_{1/\sqrt{\varepsilon} \rightarrow +\infty} d^L_2=
			- \frac{e^{- \pi\,i\,\nu}}{\Gamma(\nu)}\,.				
		\een
			
	\ele				
  
	\proof
	Let us represent $d^R_3$ as
	$$\,
	 d^R_3=-\frac{1}{2 \Gamma(\nu)} z^{1-\nu}
	\frac{\Gamma(z+\frac{\nu}{2})}{\Gamma(z)\,z^{\frac{\nu}{2}}}\,z^{\frac{\nu}{2}}\, 
	\frac{\Gamma(z)\,z^{1-\frac{\nu}{2}}}{\Gamma(z-\frac{\nu}{2}+1)}\,z^{-1+\frac{\nu}{2}}\,,
	\,$$
	where $z:=1/2 \sqrt{\varepsilon}$. Then the statement follows from the limit 
	(see \cite{HB-AE} formula 1.18(5))
	  \ben
	    \lim_{|z| \rightarrow \infty}
				\frac{\Gamma(z+\alpha)}{\Gamma(z)\,z^{\alpha}}=1\,.
		\een
		In the same manner one can prove the statement for $d^L_2$.
	\qed

	Note that for non-resonant values of the parameter $\sqrt{\varepsilon}$
	the matrices
	 \ben
	  e^{\pi\,i(\Lambda + \frac{1}{x_j} Q)}=\left(\begin{array}{ccc}
		   e^{\pm \frac{\pi\,i}{\sqrt{\varepsilon}}}   &0    &0\\[0.15ex]
			 0   &e^{\pi\,i(\nu-2\pm\frac{2}{\sqrt{\varepsilon}})}   &0\\[0.15ex]
			 0   &0          &e^{\pi\,i(\nu-4)}
			                                         \end{array}
																					\right)		
	\een
	will be divergent when $\sqrt{\varepsilon} \rightarrow 0$. But during a logarithmic resonance
	of type ${\bf (B)}$ and ${\bf (C)}$ these matrices stay constant
	 $$\,
	   e^{\pi\,i(\Lambda + \frac{1}{x_L} Q)}= 
		\diag\left(e^{-\pi\,i\,\nu},\, e^{-\pi\,i\,\nu},\,e^{\pi\,i\,\nu}\right)\,,\quad
		e^{\pi\,i(\Lambda + \frac{1}{x_R} Q)}= 
		\diag\left(e^{\pi\,i\,\nu},\, e^{3 \pi\,i\,\nu},\,e^{\pi\,i\,\nu}\right)\,,
		\,$$
		because of the relation $1/2\sqrt{\varepsilon}-\nu/2\in\NN$. So, in these cases the limits
		$\lim_{\sqrt{\varepsilon} \rightarrow 0} M_j(\varepsilon), \,j=L, R$ exist. That is why taking values of
		$\sqrt{\varepsilon}$ for which these matrices  stay constant is a good idea.

Now we can state the main result of this paper.

\bth{main}
Assume that $\nu\in\RR$ is fixed. Assume also that  
 $1/\sqrt{\varepsilon} - \nu\in 2 \NN$. Then 
		\ben
		  e^{2 \pi\,i\,T_L} \longrightarrow St_\pi\,,\qquad
			e^{2 \pi\,i\,T_R}
			  \longrightarrow St_0\,,
		\een		
			when $\sqrt{\varepsilon} \rightarrow 0$.
	 \ethe

  \proof
	
	From Theorem 4.25 of \cite{CL-CR} it follows that the unfolded Stokes matrices
	$St_L(\varepsilon)$ and $St_R(\varepsilon)$ tend to the Stokes matrices $St_{\pi}$
	and $St_0$ of the initial equation when $\sqrt{\varepsilon} \rightarrow 0$. Then from
	\prref{S-T} and the symmetry  \eqref{symm} of the perturbed equation it follows that
	the matrices $e^{2 \pi\,i\,T_L}$ and 
	$e^{2 \pi\,i\,T_R}$
  tend to the Stokes matrices $St_{\pi}$ and $St_0$ when $\sqrt{\varepsilon}$ tends to 0.
	In particular, thanks to \leref{first} we have that
	\ben
	  e^{2 \pi\,i\,T_L} \longrightarrow 
		\left(\begin{array}{ccc}
		    1    &-\frac{2 \pi\,i\,e^{-\pi\,i\,\nu}}{\Gamma(\nu)}   &0\\[0.25ex]
				0    &1                                                 &0\\[0.25ex]
				0    &0                                                 &1
				 \end{array}
		\right)=St_{\pi}		
	\een
	and
	\ben
	  e^{2 \pi\,i\,T_R}
		\longrightarrow 
		\left(\begin{array}{ccc}
		    1    &0                 &-\frac{ \pi\,i}{\Gamma(\nu)}  \\[0.25ex]
				0    &1                 &0\\[0.25ex]
				0    &0                 &1
				 \end{array}
		\right)=St_0\,,
	\een
	when $\sqrt{\varepsilon} \rightarrow 0$.
	The latter confirms one more time the statement of the theorem. This end the proof.
	\qed

																					\vspace{1cm}
			{\bf Acknowledgments.}	The author thanks the referee for valuable suggestions
			and comments, which led to the simplification and clarification of the paper.
			The author is grateful to L. Gavrilov	and E. Horozov
			for helpful  discussions and comments.
			The author was partially supported by Grant  DN 02-5/2016 
    of the Bulgarian Fond ``Scientific Research''.

\begin{small}
    
\end{small}

\end{document}